\documentclass[leqno]{article}
\usepackage{amssymb,amsfonts,amsmath,amsthm}
\usepackage[all]{xypic}
\renewcommand{\thefootnote}{\fnsymbol{footnote}} 
\newcounter{prop}
\newtheorem{prop}[subsection]{Proposition}

\newtheorem{lem}[subsection]{Lemma}

\newtheorem{cor}[subsection]{Corollary}
\newtheorem{theo}[subsection]{Theorem}
\title{The Spin of Prime Ideals}
\author{J. B. Friedlander, H. Iwaniec, B. Mazur and K. Rubin}
\date{}

\def\F{\mathbf{F}}
\def\Sel{\mathrm{Sel}}
\def\calu{\mathcal{U}}
\def\Fbb{\mathbb{F}}

\def\beq{\begin{eqnarray}}
\def\eeq{\end{eqnarray}}

\def\fa{\mathfrak a}
\def\fA{\mathfrak A}
\def\fb{\mathfrak b}
\def\fB{\mathfrak B}
\def\fc{\mathfrak c}
\def\fC{\mathfrak C}
\def\fd{\mathfrak d}
\def\fD{\mathfrak D}
\def\ff{\mathfrak f}
\def\fh{\mathfrak h}
\def\fg{\mathfrak g}
\def\fm{\mathfrak m}
\def\fM{\mathfrak M}
\def\fn{\mathfrak n}
\def\fp{\mathfrak p}
\def\fP{\mathfrak P}
\def\fq{\mathfrak q}

\renewcommand{\mod}{\mathop{\rm mod}}

\def\le{\leqslant}

\def\ge{\geqslant}

\def\ve{\varepsilon}

\def\rb{{\mathbb R}}
\def\lb{{\mathbb L}}
\def\mb{{\mathbb M}}
\def\zb{{\mathbb Z}}
\def\qb{{\mathbb Q}}
\def\sumflat{\sideset{}{^\flat}\sum}
\def\sumwedge{\sideset{}{^\wedge}\sum}

\def\ssum{\mathop{\sum\!\sum}}

\def\sssum{\mathop{\sum\!\sum\!\sum}}

\def\mychi#1{\chi_{\raise -.5ex\hbox{\the\scriptfont0 {#1}}}}

\numberwithin{equation}{section}

\begin{document}

\maketitle
\parskip=11pt
\begin{abstract}

Fixing a nontrivial automorphism of a number field $K$,  we associate 
to ideals in $K$ an invariant (with values in $\{0,\pm 1\}$) 
which we call the {\it  spin} and for which the associated $L$-function 
does not possess Euler products. We are nevertheless able, using the 
techniques of bilinear forms, to handle spin value distribution over primes, 
obtaining stronger results than the analogous ones which follow from the 
technology of $L$-functions in its current state. The initial application 
of our theorem is to the arithmetic statistics of Selmer groups of 
elliptic curves.

\end{abstract}

\renewcommand{\thefootnote}{\fnsymbol{footnote}} 
\footnotetext{\noindent{\\[-5pt]
J.\ B.\ Friedlander: Department of Mathematics, University of Toronto, M5S 2E4 Canada \\
\phantom{xxxxx}{\tt frdlndr@math.toronto.edu} \\[2pt]
H.\ Iwaniec: Department of Mathematics, Rutgers University, New Brunswick NJ 08903, USA \\
\phantom{xxxxx}{\tt iwaniec@math.rutgers.edu} \\[2pt]
B.\ Mazur: Department of Mathematics, Harvard University, Cambridge MA 02138, USA \\
\phantom{xxxxx}{\tt mazur@math.harvard.edu} \\[2pt]
K.\ Rubin: Department of Mathematics, UC Irvine, Irvine CA 92697, USA \\
\phantom{xxxxx}{\tt krubin@math.uci.edu}
\\[4pt]
{\em $2010$ Mathematics Subject Classification.} Primary 11R44; Secondary 11G05} }    
\renewcommand{\thefootnote}{\fnsymbol{footnote}}

\section{Introduction}\label{sec 1}

A very attractive area where algebraic and analytic number theory meet 
is in the distribution of prime ideals. Such questions are typically 
studied by means of $L$-functions, starting with Dirichlet and 
continued by Dedekind, Hecke, Cebotarev, Artin, and many others. The 
key feature of such $L$-functions is the Euler product, which enables analytic 
arguments, particularly analytic continuation, functional equations and 
zero-free regions to be employed. 

Given a number field $K$ and an automorphism $\sigma : K\rightarrow K$ 
we associate to ideals in $K$ an invariant for 
which the associated $L$-function does not possess Euler products; 
nevertheless we are able to handle the distribution of this 
invariant over primes. Actually, 
we obtain stronger results than the analogous ones which follow from 
$L$-function theory in its current state. We refer to our symbol as the 
``spin'' of the ideal. The spin occurs most naturally if we confine ourselves 
to number fields with specific properties. A more thorough discussion of these 
is given in Section~\ref{sec 2}. 

Let $K/\qb$ be a Galois extension of degree $n\ge 3$ with cyclic Galois group 
$G= Gal (K/\qb)$. For simplicity, we assume that $K$ is totally real and that 
the totally positive units are exactly the squares of units. At the end of 
Section~\ref{sec 2} we give some examples which show that there is 
a plentiful supply 
of fields satisfying these conditions. One convenient feature of such fields 
is the coincidence of the ideal class groups: 
class equivalence is the same, whether defined in the narrow 
or the wide sense. Although we do 
not require the class number $h$ to be one, we shall assign spins only to 
principal ideals; for these it can be done most neatly. 

For a given element $\sigma$ of $G$ and each odd principal ideal 
$\frak a$ we define its spin by 
\begin{equation}\label{eq:1.1}  
spin (\sigma, \frak a)= 
\left(\frac{\alpha}{\frak a ^{\sigma}}\right) , 
\end{equation}   
where $\frak a = (\alpha )$, $\alpha \succ 0$ (that is, $\alpha$ is 
totally positive) and $(\alpha / \frak b )$ 
stands for the quadratic residue symbol in $K$. 
For simplicity, we assume, apart from the final two sections 
of the paper, that $\sigma$ is a generator of $G$ and is fixed throughout 
so that we may denote the spin for brevity as $spin (\frak a)$, 
keeping in mind that the spin 
depends on the choice of the generator $\sigma$. Note that $spin (\frak a ) 
= \pm 1$ if $(\frak a, \frak a^{\sigma}) =1$ and is zero otherwise.  

An essential ingredient in our arguments is a bound for short sums of 
real Dirichlet characters; see Conjecture $C_n$ in Section~\ref{sec 8}. 

\begin{theo}\label{theo 1.1} Let $n\ge 3$. 
Assume Conjecture $C_n$ with exponent $\delta \le 2/n$. We have 
\begin{equation}\label{eq:1.2}  
\sum_{\substack{\frak p\,\,\, {\rm principal}\\N\frak p \le x}}
spin (\frak p) \ll x^{1-\nu +\ve} , 
\end{equation}   
where $\fp$ runs over odd prime ideals and 
$\nu =\nu (n) =\frac{\delta}{2n(12n+1)}$. 
Here, the implied constant depends on 
$\ve$ and the field $K$.
\end{theo}

It turns out that our Conjecture $C_3$ holds true with exponent 
$\delta (3)=1/48$ due to a well-known result of D. Burgess [Bu]. 
Hence, the estimate~\eqref{eq:1.2} is unconditionally true for cubic 
fields with  
\begin{equation}\label{eq:1.3}  
\nu = \nu (3) = 1/10656 .  
\end{equation}   

The contribution to the sum in~\eqref{eq:1.2} coming from primes 
of degree greater than one is negligible and, for unramified primes of 
degree one, the spin takes the values $\pm 1$. Our theorem implies that, 
unconditionally for $n=3$ and, dependent on Conjecture $C_n$ for $n\ge 4$, 
the spin takes these values asymptotically equally often. 
Hopefully, future progress with short character sums will 
make~\eqref{eq:1.2} unconditional for fields of higher degree.

For $n = 2$ the sign change of the spin of principal prime ideals of 
degree one is quite regular, so our method (for catching primes by  
bilinear forms techniques) fails (however, see further results on this  
in the final Section~\ref{sec 10}). 
In fact, if we restrict the primes to a fixed residue class of 
modulus $8\Delta$ where $\Delta$ is the field discriminant, then the spin 
is constant, so for such a restricted sum~\eqref{eq:1.2} is false.  
On the other hand, if $n \ge 3$ and $\sigma$ is a generator,  
then the spin turns quite randomly, so~\eqref{eq:1.2} can be 
established for sums twisted by many kinds of characters. For example, 
our arguments work with very little change when the spin is twisted by 
a Hecke Gr\"ossencharakter.

Note that the theorem saves a fixed power of $x$ when summing the prime spins. 
This is in contrast to the $L$-function theory which, in the absence of 
something approaching the Riemann Hypothesis, would not permit us to even 
count these primes with such a degree of accuracy. This observation suggests 
that the sets of primes of constant spin, although determined by a natural 
algebraic condition, seem very unlikely to be 
{\it Cebotarev classes} (that is, sets of primes distinguished by 
splitting properties in some fixed finite extension of $K$). 

As it happens, there is little extra work involved in proving 
a more general result showing that cancellation occurs 
when summing the spin over primes in an arithmetic progression and 
this is of interest both on its own and for applications. 
Let $\fM$ be an integral ideal of $K$ with $2\mid \fM$ and $\mu$ an 
integer of $K$ with $(\mu, \fM)=1$. We shall take the progression 
$\mu\, (\mod\, \fM)$ to be fixed throughout the paper. 

\begin{theo}\label{theo 1.2}
The bound~\eqref{eq:1.2} still holds when the sum is further restricted 
to principal prime ideals which have a totally positive generator $\pi$ 
satisfying $ \pi \equiv \mu\, (\mod \fM)$.   
\end{theo}

Our results have implications for the distribution of the Selmer rank of 
elliptic curves. In fact, the question addressed in this paper first arose 
in trying to improve the results in [MR] on Selmer ranks in families of 
quadratic twists. 
Let $E$ for example
 be the elliptic curve 
$$y^2 = x^3 + x^2 -16x-29\ ,
$$
which has conductor $784 = 2^4\cdot 7^2$.  
Let $K$ be the maximal real subfield of the field $\qb(\boldsymbol{\mu}_7)$ 
of $7$-th roots of unity.  
Then, $K$ is a cyclic extension of $\qb$ of degree $3$, and $K = \qb(E[2])$, 
the field generated by the coordinates of the points of order $2$ on $E$.

Suppose $p$ is a rational prime congruent to $\pm 1 \pmod{7}$, so 
$p$ splits into $3$ distinct primes in $K$.  
Let $\fp$ be one of the primes above $p$.
If $\fp$ has a totally positive generator that is congruent to $1 \pmod{8}$, 
then the $2$-Selmer group $\Sel_2(E^{(p)}/\qb)$ of the quadratic twist of $E$ 
by $p$ has dimension
$$
\dim_{\F_2}\Sel_2(E^{(p)}/\qb) = 
\begin{cases}
3 & \text{if $spin(\fp) = \, 1$}, \\
1 & \text{if $spin(\fp) = -1$}.
\end{cases}
$$
The condition that $p$ have a generator congruent to 1 modulo 8 is 
equivalent to asking that $p$ split completely in the ray class field 
of $K$ modulo 8. Hence, the set of such $p$ has positive density. 
Moreover, $K$ has class number $1$. Thus,  
Theorem~\ref{theo 1.2} shows that, within that set of twists, 
the Selmer rank is equal to 1 half of the time and 
3 half of the time. As one might expect, this holds more generally; 
see Section~\ref{sec 11}. 

We conclude the introduction with a brief outline of the contents of the 
paper. In Section~\ref{sec 2.0} we recall the law of quadratic reciprocity 
in the setting of a general number field and some related issues which 
we need for this work. In Section~\ref{sec 2} we recall some basic 
facts about number fields; to some extent these are specialized to the 
type of fields we are considering. One 
of the technical problems we encounter is the difficulty in stepping 
smoothly between ideals and integers. This requires a good understanding 
of the geometrical shape of a convenient fundamental domain for the action 
of units on $\rb^n$. In Section~\ref{sec 3} we use the construction of 
T. Shintani [Shi] and we establish various properties which are needed 
for our applications. 

Although the quadratic residue symbol has multiplicative properties, the 
spin does not behave in a purely multiplicative fashion. This feature 
rules out the possibility of using the techniques of $L$-functions but opens 
the possibility of using the well-known technique of transforming sums 
over primes to congruence sums and bilinear forms. In Section~\ref{sec 4} 
we show how sums over primes are reduced to these latter shapes in a 
general context and then, in 
Sections~\ref{sec 5} and~\ref{sec 6} respectively, we produce the 
required bounds for them.  
Section~\ref{sec 7} quickly pieces together these ingredients to complete the 
proof of the theorems. In Section~\ref{sec 8} 
we present the estimates for character sums which power the whole work. 
In particular, for cubic fields ($n=3$), the bound we use was established 
by Burgess via an appeal to the Riemann Hypothesis for algebraic curves,  
which in this case have genus six. 

In Section~\ref{sec 11} we prove Theorem~\ref{thm11.1} which 
relates the spin to Selmer groups and justifies, in greater generality, 
the claims made for our example. 
In Section~\ref{sec 9} we briefly discuss the relationships 
among the spins associated to different automorphisms and 
point out a few of the very 
interesting (to us) problems left open by this work.  Then, in the 
final Section~\ref{sec 10} we consider prime spins for an involutary 
automorphism, a problem which requires completely different tools.

{\bf Acknowledgements}: This work was initiated and a significant part 
was accomplished  
during the time the authors enjoyed participating, with MSRI 
financial support, in the program on  
``Arithmetic Statistics'', held at MSRI Berkeley during January--May,  
2011. The second and third-named authors also received support from the 
Clay Mathematics Institute. Research of J. F. is supported by 
NSERC grant A5123, 
that of H. I. by NSF Grant DMS-1101574, that of B. M. 
by NSF Grant DMS-0968831 and that of K. R. by NSF Grant 
DMS-1065904. 
 We would also like to thank the referee for a very thorough reading 
of the paper. 

\section{Quadratic Residues and Reciprocity}\label{sec 2.0}

We say that an integral ideal $\frak a$ is odd if $(\frak a, 2)=1$,  
and is otherwise even. 
If $\frak p$ is an odd prime ideal and $\alpha$ is an integer in $K$ with 
$\alpha\not \equiv  0\, (\, \mod \frak p)$, then the quadratic residue 
symbol $(\alpha/\frak p)$ is defined by 
$$
\left(\frac {\alpha}{\frak p}\right)\,  =  \, 
\begin{cases}\,\,\, 1 
 \quad \text{if $\alpha \equiv \xi^2 \, (\mod \frak p)$}\, , 
\\
{-1} 
 \quad \text{otherwise}\, .
\end{cases}
$$
We extend this to all $\alpha \in \cal O$, the ring of integers of $K$, 
by setting 
$$\left(\frac {\alpha}{\frak p}\right)\,  =  0 \,  
\quad \text{if $\alpha\equiv  0\, (\, \mod \frak p)$}\, .
$$
If $\frak q = \frak p_1\ldots\frak p_r$ is the product of odd prime ideals 
(not necessarily distinct) then 
$$
\left(\frac {\alpha}{\frak q}\right)\,= \left(\frac {\alpha}{\frak p_1}\right)
\ldots\left(\frac {\alpha}{\frak p_r}\right)
$$
gives us a multiplicative extension of the definition. 
Given an odd ideal $\frak q$, the symbol $(\alpha/\frak q)$ is periodic 
in $\alpha \, (\mod \, \frak q)$ and multiplicative in $\alpha$. 

Have in mind that the quadratic residue symbol $(\alpha/\frak q)$ depends 
on the field $K$ although this fact is not displayed in our notation. 
Observe that, for $p$ a rational prime not dividing the rational integer $a$, 
the symbol $(a/p)$ does not agree with the Legendre symbol in $\qb$. To see 
this, take a prime $p$ which splits completely in $K$. Suppose $K/\qb$ is 
Galois of even degree. Then we can write $p=\frak q\frak q^{\sigma}$ where 
$\sigma$ is an involution in $G=Gal(K/\qb)$. Hence, 
$$
\Bigl(\frac{a}{p}\Bigr) =\Bigl(\frac{a}{\frak q}\Bigr)\,\,
\Bigl(\frac{a}{\frak q^{\sigma}}\Bigr) = \Bigl(\frac{a}{\frak q}\Bigr)^2 
= 1\ , 
$$
while the Legendre symbol sometimes changes values. 

If 
$\beta \in \cal O$, we say that $\beta$ is odd or even according as 
$(\beta)$ is, where $(\beta)$ is the prinicipal ideal generated by $\beta$, 
and for $\beta$ odd, we define the symbol $(\alpha/\beta)  = 
(\alpha/(\beta))$.
Note that $(\alpha/\beta)$ does not change if $\beta$ is replaced by an 
associate; however, it may vary if $\alpha$ is changed by a unit which 
is not a square. 

The famous law of quadratic reciprocity extends to arbitrary number fields. 
A nice treatment of this is given in Chapter 8 of Hecke [He]. See also 
Chapter 6, Section 8 of J. Neukirch [Ne]. To state this in the most 
convenient way, 
we first require the Hilbert symbols. 

For any $\alpha,\,\beta\in K$ and any $\frak p$, the Hilbert symbol 
$$
\Bigl(\frac{\alpha,\,\beta}{\frak p}\Bigr) = \pm 1\ ,  
$$
takes the value $1$ if the quadratic form 
$\alpha x^2 +\beta y^2 -z^2$ represents zero non-trivially in the local 
field $K_{\fp}$ and takes the value $-1$ otherwise. 
We have the following properties for any 
$\alpha,\,\beta,\, \gamma\in K$ and any place $\fp$ (see Proposition 3.2 
of [Ne]): 
\begin{equation*}
\begin{split}
\Bigl(\frac{\alpha\, ,\beta}{\fp}\Bigr)\,  & = 
\Bigl(\frac{\beta\, ,\alpha}{\fp}\Bigr)\ ,\quad 
\Bigl(\frac{\alpha\, ,\beta\gamma}{\fp}\Bigr) = \,  
\Bigl(\frac{\alpha\, ,\beta}{\fp}\Bigr)\, 
\Bigl(\frac{\alpha\, ,\gamma}{\fp}\Bigr)\ , \\
\Bigl(\frac{\alpha\, , -\alpha}{\fp}\Bigr) & = 1\, , \quad  
\Bigl(\frac{\alpha\, , 1 -\alpha}{\fp}\Bigr) = 1\, , \quad 
\Bigl(\frac{\alpha\, , \beta^2}{\fp}\Bigr) = 1\, .
\end{split}
\end{equation*}

We next define, for any  
$\alpha,\,\beta,\, \in \cal O$ 
$$
\mu_2(\alpha,\beta) = 
\prod_{\frak p \mid 2}\Bigl(\frac{\alpha,\,\beta}{\frak p}\Bigr)\, ,\quad
\mu_{\infty}(\alpha,\beta)= 
\prod_{\frak p \mid \infty}\Bigl(\frac{\alpha,\,\beta}{\frak p}\Bigr)\ ,
$$
and
$$
\mu(\alpha,\beta) =  \mu_2(\alpha,\beta)\, \mu_{\infty}(\alpha,\beta)\ .
$$

For any $\alpha,\, \beta \in {\cal O}$, $\beta$ odd, we introduce 
the {\it completed} quadratic residue symbol 
\begin{equation}\label{eq:2.5.1}
\Big|\frac{\alpha}{\beta}\Big| = 
\mu_{\infty}(\alpha,\beta)\, \Bigl(\frac{\alpha}{\beta}\Bigr) \ . 
\end{equation}
Note that, if $\alpha$ or $\beta$ is totally positive, then
$|\alpha/\beta| = (\alpha/\beta)$. 

Now, we can state: 
\begin{lem}\label{lem 2.0}(Law of Quadratic Reciprocity) 
For two odd integers $\alpha$, $\beta\, \in \cal O$, we have
\begin{equation}\label{eq:2.4}
\Bigl(\frac {\alpha}{\beta}\Bigr)\,  = \mu(\alpha,\beta)\, 
\Bigl(\frac {\beta}{\alpha}\Bigr)\ ,
\end{equation}
or equivalently, 
\begin{equation}\label{eq:2.5}
\Bigl|\frac {\alpha}{\beta}\Bigr|\,  = \mu(\alpha,\beta)\, 
\Bigl|\frac {\beta}{\alpha}\Bigr| 
= \mu_2(\alpha,\beta)\, 
\Bigl(\frac {\beta}{\alpha}\Bigr)\ . 
\end{equation}
\end{lem}

Even in the rational field, the reciprocity law finds problems with
the prime $2$ and in number fields that situation is further complicated. 
We shall need to circumvent some of those problems with the symbol 
$(\alpha/\beta)$ in the situation where the upper entry may not be 
odd and we could not find in the literature a treatment which completely 
fulfilled our needs. 

We shall show 
that the symbol $|\alpha/\beta|$, as a function of $\beta$ is periodic 
of period $(8\alpha)$.
Actually, this does what we require   
with room to spare; all we need is that, for $(\alpha,\,\beta)=1$ 
the symbol $(\alpha/\beta)=\pm 1$ depends on $\alpha$ but only on the 
residue class of $\beta$ modulo $(2^{\ell}\alpha)$, for some $\ell$ 
which could even depend on the field $K$. See the use of this property 
in Section~\ref{sec 5}. 

Our goal in the remainder of this section is: 
\begin{prop}\label{cor 2.1.2} 
Fix any nonzero $\alpha \in \cal O$. The symbol $|\alpha/\beta|$, for 
odd integers $\beta$, depends only on the residue class of $\beta$ 
modulo $(8\alpha)$. The same is true for the symbol $(\alpha/\beta)$ 
for odd totally positive integers $\beta$. 
\end{prop} 

We begin the argument with the following result. 
\begin{lem}\label{lem 2.1.1}
Let $\alpha,\,\beta\in {\cal O}$ be odd. 
Then, $\mu_2(\alpha, \beta)$ depends only on the residue classes of 
$\alpha, \, \beta$ modulo $8$. 
\end{lem}
\begin{proof}
Using the multiplicativity property of the Hilbert symbol, we see 
that it suffices to prove that any $\eta \equiv 1 (\mod 8)$ 
is a square in $K_{\fp}$ for any $\fp | 2$. This  
follows by Hensel's lemma, but actually, it is seen explicitly   
from the identity 
$$\sqrt{1+8x} = 1 + 
\sum_{\ell \ge 1}c_{\ell}
\frac{(4x)^{\ell}}{\ell \, !}\ , \quad \rm {where} \,\,c_\ell 
= \prod_{0\le k <\ell}(1-2k)\ .
$$ 
Since $2^{\ell}$ does not 
divide $\ell \, ! $, the series converges $\fp$-adically.  
This completes the proof of the lemma. 
\end{proof}
Combining the lemma with the reciprocity law, we obtain 
\begin{cor}\label{cor 2.1.3} 
Proposition~\ref{cor 2.1.2} is true in case $\alpha$ is odd. 
\end{cor}

So, as expected, the main problem occurs when $\alpha$ is even. 
Were we to have a reasonable definition for our symbol when the 
lower entry is even and an accompanying version of the reciprocity 
law, then this case would probably also 
be straight-forward. As it is, we manoeuvre to 
reduce to the situation of the upper entry being odd by using properties of 
the symbol, in particular, its already known periodicity.

We next consider a further supplement to the Legendre symbol for which 
the period is even smaller. We define, for $\alpha, \, \beta \in \cal O$, 
 $\beta$ odd,  
\begin{equation}\label{eq:2.5.3}
\Big[\frac{\alpha}{\beta}\Big] = \mu_2(\alpha,\,\beta)
\Bigl|\frac{\alpha}{\beta}\Bigr| = \mu(\alpha,\,\beta)
\Bigl(\frac{\alpha}{\beta}\Bigr) \ .
\end{equation}

\begin{lem}\label{lem 2.1.2} 
Fix $\alpha \in \cal O$ such that $1+\alpha$ is odd. Then 
\begin{equation}\label{eq:2.5.4}
\Big[\frac{\alpha}{\beta}\Big] = \Big[\frac{\alpha}{\delta}\Big] 
\quad {\rm if}\,\, \beta\equiv \delta\, (\mod (2\alpha))\ .
\end{equation}
\end{lem}

\begin{proof}
We can assume $(\alpha\, , \beta\delta)=1$ since 
otherwise~\eqref{eq:2.5.4} trivially holds.  Fix $\gamma \in {\cal O}$ 
such that $\beta\,\gamma\equiv 1+\alpha\, (\mod (2\alpha))$. 
Note that $\gamma$ is odd. It suffices to show that
$$\Big[\frac{\alpha}{\beta}\Big] = \Big[\frac{\alpha}{\gamma}\Big] \ .
$$
To this end, consider the number 
$$
\lambda=\frac{\beta\,\gamma -1}{\alpha}\, , \quad {\rm so}\,\,\, 
\beta\,\gamma =\alpha\,\lambda + 1\ .
$$
Obviously $\lambda \in {\cal O}, \,\, \lambda \equiv 1\, (\mod\, 2)$ 
and $(\lambda\, , \, \beta\,\gamma )=1$. We write 
\begin{equation*}
\begin{split}
\Bigl(\frac{\alpha}{\beta}\Bigr)\, =
\Bigl(\frac{\alpha}{\beta}\Bigr)\,\Bigl(\frac{\lambda}{\beta}\Bigr)^2\,  
& = \Bigl(\frac{\alpha\lambda }{\beta}\Bigr)\,\Bigl(\frac{\lambda}{\beta}
\Bigr)\, = \Bigl(\frac{-1}{\beta}\Bigr)\,\Bigl(\frac{\lambda}{\beta}
\Bigr)\, = \Bigl(\frac{- \lambda}{\beta}\Bigr)\\  
& = \Bigl(\frac{\beta}{\lambda}\Bigr)
\prod_{\frak p \mid 2\infty}\Bigl(\frac{\beta\, ,\,-\lambda}{\frak p}\Bigr)
\ ,
\end{split}
\end{equation*}
by the reciprocity law. Next, we write 
$$
\Bigl(\frac{\beta}{\lambda}\Bigr)= \Bigl(\frac{\gamma}{\lambda}\Bigr)
= \Bigl(\frac{\lambda}{\gamma}\Bigr)
\prod_{\frak p \mid 2\infty}\Bigl(\frac{\gamma,\,\lambda}{\frak p}\Bigr)\ . 
$$ 
Here, we have
$$
\Bigl(\frac{\lambda}{\gamma}\Bigr)
=\Bigl(\frac{\alpha}{\gamma}\Bigr)\,\Bigl(\frac{\alpha\lambda}{\gamma}\Bigr)
=\Bigl(\frac{\alpha}{\gamma}\Bigr)\,\Bigl(\frac{-1}{\gamma}\Bigr)
=\Bigl(\frac{\alpha}{\gamma}\Bigr)\,
\prod_{\frak p \mid 2\infty}\Bigl(\frac{\gamma,\,- 1}{\frak p}\Bigr)\ .
$$
Collecting the above results, we arrive at 
$$
\Bigl(\frac{\alpha}{\beta}\Bigr)\, =\Bigl(\frac{\alpha}{\gamma}\Bigr)\, 
\prod_{\frak p \mid 2\infty}\mu(\fp)\ ,
$$
with 
\begin{equation*}
\begin{split}
\mu(\fp) & = 
\Bigl(\frac{\beta\, ,\,-\lambda}{\frak p}\Bigr)\,
\Bigl(\frac{\gamma\, ,\,\lambda}{\frak p}\Bigr)\, 
\Bigl(\frac{\gamma\, ,\,- 1}{\frak p}\Bigr)\, 
= \Bigl(\frac{\beta\gamma\, ,\,-\lambda}{\frak p}\Bigr) \\
& = \Bigl(\frac{\beta\gamma\, ,\,-\alpha\lambda}{\frak p}\Bigr)
\Bigl(\frac{\beta\gamma\, ,\,\alpha}{\frak p}\Bigr) 
= \Bigl(\frac{\beta\gamma\, ,\,1-\beta\gamma}{\frak p}\Bigr)
\Bigl(\frac{\beta\gamma\, ,\,\alpha}{\frak p}\Bigr) \\
& = \Bigl(\frac{\beta\gamma\, ,\,\alpha}{\frak p}\Bigr)
= \Bigl(\frac{\beta\, ,\,\alpha}{\frak p}\Bigr)\,
\Bigl(\frac{\gamma\, ,\,\alpha}{\frak p}\Bigr)\ .
\end{split}
\end{equation*}
This completes the proof of~\eqref{eq:2.5.4} 
and the lemma.
\end{proof} 

We are now ready to complete the proof of Proposition~\ref{cor 2.1.2}. 
As particular cases of the previous lemma we get $[2/\beta] =[2/\delta]$ if 
$\beta\equiv \delta\, (\mod 4)$ and, for any $\alpha \neq 0$, 
$[2\alpha/\beta] =[2\alpha/\delta]$ if 
$\beta\equiv \delta\, (\mod (4\alpha))$. Multiplying the last two 
equations we find (because $4$ is a square)
$$
\Big[\frac{\alpha}{\beta}\Big] =\Big[\frac{\alpha}{\delta}\Big] 
\quad {\rm if}\,\,\, 
\beta\equiv \delta\, (\mod (4\alpha))\ .
$$
Here, if we have the stronger congruence condition 
$\beta\equiv \delta\, (\mod (8\alpha)),$ 
the Hilbert symbol at any $\fp | 2$ can be omitted 
by Lemma~\ref{lem 2.1.1}.  
This completes the proof of the proposition.

\section{Number Field Preliminaries}\label{sec 2}

Let $K$ be a totally real number field of degree $n=[K:\qb]$ so $K$ has 
$n$ embeddings into 
$\rb$. For $\alpha \in K$ we denote its conjugates by
$\alpha^{(1)},\ldots,\alpha^{(n)}$. They are all real and
$$
N\alpha = \alpha^{(1)}\ldots\alpha^{(n)}, \quad T\alpha = 
\alpha^{(1)} +\ldots +\alpha^{(n)}
$$ are the norm and the trace of $\alpha$. We say that $\alpha\in K$ 
is totally positive if all its conjugates are positive, in which case 
we write $\alpha \succ 0$. We embed $K$ into $\rb^n$ by the mapping 
\begin{equation}\label{eq:2.1}
\alpha \rightarrow (\alpha^{(1)},\ldots,\alpha^{(n)})
\end{equation} 
with addition and multiplication performed componentwise. With a slight 
abuse of notation, for ${\cal D}$ a subset of $\rb^n$ we shall write 
briefly $\alpha \in {\cal D}$ meaning 
$(\alpha^{(1)},\ldots,\alpha^{(n)}) \in {\cal D}$. 

Let $\mathcal U$ denote the group of units, $\mathcal U^+$ the subgroup of 
totally positive units and $\mathcal U^2$ the subgroup of 
squares of units, so $\mathcal U^2\subset \mathcal U^+ \subset 
\mathcal U$. Note that $[{\cal U}: {\cal U}^2] = 2^n$ because, 
by the Dirichlet unit theorem, every $u\in \cal U$ has a unique 
representation 
\begin{equation}\label{eq:2.2}
u = \pm\, \ve_1^{m_1}\ldots\ve_r^{m_r}, \,\,\, r=n-1, 
\end{equation} 
where $\ve_1, \ldots, \ve_r$ is a system of fixed fundamental units 
and $m_1, \ldots , m_r \in \zb$. Hence $u$ is a square 
exactly when it has positive sign and the exponents $m_1, \ldots , m_r$ 
are even. 
Assume that the homomorphism $\mathcal U\rightarrow \{\pm 1\} \times 
\ldots \times \{\pm 1\}$ given by 
$$
u \rightarrow 
\left(\frac{u^{(1)}}{ |u^{(1)}|}, \ldots , \frac{u^{(n)}}{ |u^{(n)}|}\right)
$$
is surjective. 
Equivalently, this means that $[{\cal U}: {\cal U}^+] = 2^n$, hence 
${\cal U}^+ = \mathcal U^2$ and also that ideal class equivalence is the 
same, whether defined in the wide or the narrow sense. 

Note that for $u\in {\cal U}^+,\, u\neq 1$ we have 
\begin{equation}\label{eq:2.3}
u^{(1)} + \ldots + u^{(n)} > n 
\end{equation} 
and equality holds for $u=1$. This follows from the 
arithmetic-geometric mean inequality. 

We are now going to define the spin of odd principal ideals. Throughout, 
we assume that $K/\qb$ is a totally real Galois cyclic extension of 
degree $n\ge 3$ and ${\cal U}^+ ={\cal U}^2$. Fix a generator of 
$G= Gal (K/\qb)$, say $\sigma$. Then, if $\frak a$ is an odd principal 
ideal, we define 
\begin{equation}\label{eq:2.7}
spin (\frak a) = \Bigl(\frac{\alpha}{\frak a^{\sigma}}\Bigr)\ ,
\end{equation}
where $\alpha$ is chosen as a totally positive generator of $\frak a$. 
Such an $\alpha$ is uniquely determined up to the square of a unit so 
$spin (\frak a)$ is well-defined. \footnote{For Gaussian primes, 
the name ``spin'' was used in [FI], but for a 
symbol which is only superficially reminiscent of our spin$(\fp)$ 
for prime ideals. 
Writing $\pi =r + is \in \zb[i]$ uniquely with $r,\, s >0, \, r$ odd, 
the spin of $p=\pi\bar\pi$ was defined to be the (usual) Jacobi symbol 
$\sigma_p =(s/r)= \pm 1$.}

Although the definition~\eqref{eq:2.7} makes sense for any $\sigma \in G$,  
we decided to choose $\sigma$ from the generators of $G$ because 
it will be convenient for our arguments. Having fixed $\sigma$, for 
notational convenience we shall write, for any ideal $\fa$ 
$$
\fa^{\sigma} =\fa', \quad \fa^{\sigma^{-1}} =\fa^- \ ,
$$
and for any $\alpha \in K$, 
$$
\alpha^{\sigma} =\alpha', \quad \alpha^{\sigma^{-1}} =\alpha^- \ .
$$
In this notation~\eqref{eq:2.7} becomes
$$
spin (\frak a) = \Bigl(\frac{\alpha}{\alpha'}\Bigr)
= \Bigl(\frac{\alpha^-}{\alpha}\Bigr)\ . 
$$
Note that $spin(\fa)= \pm 1$ if $(\fa,\fa')=1$ and   
$spin(\fa)= 0$ otherwise. 

We shall need to understand the spin of the product of two odd ideals 
which may not be principal even though their product is. To this end 
we fix a collection ${\cal C}\ell= \{\fA,\, \fB, \ldots\}$ of odd ideals, 
a set of representatives of the ideal class group, 
choosing two from each class. Put 
\begin{equation}\label{eq:2.8}
\ff = \prod_{\fC\in {\cal C}\ell}\fC\ .  
\end{equation}
We can assume, for purely technical convenience, that $\ff$ is a 
square-free ideal and of course $\ff$ is principal. Actually, 
we can assume even more, that 
\begin{equation}\label{eq:2.9}
f= N\ff \quad {\rm is\,\,\, square-free}.  
\end{equation} 
Note that this implies that the ideals in our collection, 
together with their conjugates, are pairwise co-prime and odd.
The reason for taking, in the collection ${\cal C}\ell$,  
two representatives from each ideal class is that sometimes we pick up 
one representative and find we need to do it again. For convenience, 
it is nice to have the second choice co-prime with the first. 

Now, let $\fa\fb$ be a principal ideal co-prime with $2f$. We have
\begin{equation}\label{eq:2.9.1}
\begin{split}
& \fa\fA  = (\alpha), \quad \alpha\succ 0\ , \\
& \fb\fB  = (\beta), \quad \beta\succ 0\ ,
\end{split}
\end{equation}
for some $\fA,\,\fB \in {\cal C}\ell, \, \fA\neq\fB$. Note that 
$(\fA,\,\fB)=1$ and $\fA\,\fB$ is co-prime with $\fA'\,\fB'$. 
Since $\fa\fb$ is principal so is
$\fA\fB$, say $\fA\fB = (\gamma),\,\, \gamma \succ 0$. Then 
$\fa\fb(\gamma)=(\alpha\beta)$ and 
\begin{equation*}
\begin{split}
spin (\fa\fb) &  = \Bigl(\frac{\alpha\beta\gamma}{\fa'\fb'}\Bigr)
= \Bigl(\frac{\alpha}{\fb'}\Bigr)
\Bigl(\frac{\beta}{\fa'}\Bigr)
\Bigl(\frac{\alpha}{\fa'}\Bigr)
\Bigl(\frac{\beta}{\fb'}\Bigr)
\Bigl(\frac{\gamma}{\fa'\fb'}\Bigr) \\
& = \Bigl(\frac{\alpha}{\beta'}\Bigr)
\Bigl(\frac{\beta}{\alpha'}\Bigr)
\Bigl(\frac{\alpha}{\fa'\fB'}\Bigr)
\Bigl(\frac{\beta}{\fb'\fA'}\Bigr)
\Bigl(\frac{\gamma}{\fa'\fb'}\Bigr) \\
& = \Bigl(\frac{\alpha}{\beta'}\Bigr) 
\Bigl(\frac{\beta}{\alpha'}\Bigr)
\Bigl(\frac{\alpha\gamma}{\fa'\fB'}\Bigr)
\Bigl(\frac{\beta\gamma}{\fb'\fA'}\Bigr)
\Bigl(\frac{\gamma}{\gamma'}\Bigr)\ .
\end{split}
\end{equation*}
Next, by the reciprocity law we write
$$
\Bigl(\frac{\beta}{\alpha'}\Bigr) =\Bigl(\frac{\beta^-}{\alpha}\Bigr) 
=\pm \Bigl(\frac{\alpha}{\beta^-}\Bigr)
$$
where the sign depends only on the residue classes of $\alpha,\,\beta$ 
modulo 8. Hence, we conclude the following factorization rule for the spin:
\begin{equation}\label{eq:2.10}
spin(\fa\, \fb) =\pm \Bigl(\frac{\alpha}{\beta'\beta^-}\Bigr)
\Bigl(\frac{\alpha\gamma}{\fa'\fB'}\Bigr)
\Bigl(\frac{\beta\gamma}{\fb'\fA'}\Bigr)
spin(\gamma)\ .
\end{equation}
The two middle symbols separate $\fa$ from $\fb$ and hence do not play 
a role in the estimation of general bilinear forms. However, the leading 
symbol
\begin{equation}\label{eq:2.11}
\Bigl(\frac{\alpha}{\beta'\beta^-}\Bigr)
= 
\Bigl(\frac{\alpha}{\fb'\fb^-}\Bigr)
\Bigl(\frac{\alpha}{\fB'\fB^-}\Bigr)
\end{equation}
is vital for this. Note that, if we had chosen $\sigma$ in~\eqref{eq:2.7} 
to be an involution, then $\beta'=\beta^-$ and the leading symbol 
in~\eqref{eq:2.11} would be constant, in other words the spin symbol would 
be essentially multiplicative. In such a case there is no way to get 
cancellation in general bilinear forms. This is the reason why our 
arguments fail for quadratic fields. 

If $\fp$ is an odd prime ideal then
$$
\sum_{\alpha\,(\mod \fp)}\Bigl(\frac{\alpha}{\fp}\Bigr) = 0\ .
$$
which means that the number of quadratic residue classes $(\mod \fp)$
 is equal to the number of quadratic non-residue classes $(\mod \fp)$. 
More generally, if the odd ideal $\fq$ is not the square of an ideal, 
then  
\begin{equation}\label{eq:2.12}
\sum_{\alpha\,(\mod \fq)}\Bigl(\frac{\alpha}{\fq}\Bigr) = 0\ .
\end{equation}

We say that the positive rational integer $q$ is squarefull if 
$p|q $ implies $p^2|q$. 
\begin{lem}\label{lem 2.1}
Let $\fq$ be an odd ideal whose norm $q=N\fq$ is not squarefull.
Then
\begin{equation}\label{eq:2.13}
\sum_{\alpha\,(\mod \fq'\fq^-)}\Bigl(\frac{\alpha}{\fq'\fq^-}\Bigr) = 0\ .
\end{equation}
\end{lem}
\begin{proof}
Let $p$ be a prime divisor of $q$ whose square does not divide $q$. 
Hence, $p = N\fp$ for some $\fp | \fq$ and $\fq=\fp\fc$ with $(\fc,p)=1$. 
In particular, $(\fc'\fc^-, \fp'\fp^-)=1$ so 
$\fq'\fq^- =\fp'\fp^-\fc'\fc^-$ is not the square of an ideal because 
$\fp'\neq \fp^-$. Hence~\eqref{eq:2.13} follows from~\eqref{eq:2.12}.
\end{proof}

In many situations we shall employ an integral basis of $K$, say 
$\omega_1,\ldots,\omega_n$. For convenience we can take $\omega_1=1$, so 
$$
{\cal O}= \omega_1\zb + \ldots +\omega_n\zb =\zb + \mb
$$
where $\mb = \omega_2\zb + \ldots +\omega_n\zb$ is a submodule of $\cal O$ 
of rank $n-1$. We shall also consider the submodule 
$\lb = \eta_2\zb + \ldots +\eta_n\zb$ with
\begin{equation}\label{eq:2.13.1}
\eta_2 =\omega_2-\omega'_2, \ldots, \eta_n =\omega_n-\omega'_n\ .
\end{equation}

\begin{lem}\label{lem 2.4}
The map $\mb\rightarrow {\cal O}$ given by $\beta\rightarrow\beta - \beta'$ 
is an injection; its image is the module
$\lb$. 
\end{lem}
\begin{proof}
Since the map $\xi\rightarrow \xi'$ generates the Galois group, it 
follows that, if $\beta' -\beta =0$ then $\beta$ must be rational, hence 
contained in $\zb \cap \mb$. But, $\zb \cap \mb = 0$.
\end{proof} 
\begin{cor}\label{cor 2.5}
The numbers $\eta_2,\ldots,\eta_n$  
are linearly independent over $\qb$, so $\lb$ has rank $n-1$.  
\end{cor}

It would be nice to have an integral basis
$1, \omega_2, \ldots, \omega_n$ of $K$ for which 
all the conjugates of $\theta =\eta_2/\eta_3$ are distinct, so 
\begin{equation}\label{eq:2.14}
W=\prod_{\tau \neq {\rm id}}
\bigl(\eta_2\eta_3^{\tau} - \eta_2^{\tau} \eta_3\bigr) \neq 0\ .
\end{equation}
In other words, we wish to have a basis such that $K =\qb (\theta)$. 
If $n$ is a prime number, then any basis gives $W\neq 0$. 
Indeed, if $\theta =\theta^{\tau}$ for some 
$\tau \neq {\rm id}$ then all the conjugates of $\theta$ are equal, because any
$\tau \neq {\rm id}$ is a generator of the whole Galois group $G$. Therefore, 
$\theta$ is rational, contradicting the fact that $\eta_2,\, \eta_3$ are 
linearly independent over $\qb$. However, we can work with any basis and 
any $n\ge 3$ due to the following properties.  

\begin{lem}\label{lem 2.5.1}
Let $\tau\in G\, , \, \tau \neq \,\, {\rm id}$. Then, at least two 
of the numbers
\begin{equation}\label{eq:2.14.1}
\eta_k^{\tau}/\eta_k, \quad 2\le k\le n\ ,
\end{equation} 
are distinct.
\end{lem}
\begin{proof}
Suppose all the 
numbers~\eqref{eq:2.14.1} are equal, so $(\eta_k/\eta_2)^{\tau}
=\eta_k/\eta_2 $ for $2\le k \le n$. This shows that all the $n-1$ numbers 
$\eta_k/\eta_2 $ are in a proper subfield of $K$ of degree $\le n/2 < n-1$ so 
these numbers must be linearly dependent over $\qb$, which contradicts 
Corollary~\ref{cor 2.5}. 
\end{proof}

Let $\fh(\tau)$ be the ideal generated by the numbers
\begin{equation}\label{eq:2.14.6}
\eta_k^{\tau}\eta_{\ell}-\eta_k\eta_{\ell}^{\tau}\ , 
\quad 2\le k\neq\ell \le n\ . 
\end{equation} 
Thus, Lemma~\ref{lem 2.5.1} says that $\fh(\tau)$ is a non-zero ideal. 
We denote 
\begin{equation}\label{eq:2.14.7}
\fh = \fD\prod_{\tau\ne\, {\rm id}}\fh(\tau)\ ,
\end{equation} 
where $\fD$ is the {\it different} of the field $K$.

For $\alpha = a_1\,\omega_1 + \ldots + a_n\,\omega_n \in {\cal O}$ the basis 
coefficients $a_1, \ldots, a_n$ are linear combinations of the conjugates 
$\alpha^{(1)}, \ldots,\alpha^{(n)}$ and vice-versa. Therefore, the system of 
estimates $\alpha^{(1)}, \ldots,\alpha^{(n)}\ll y$ is equivalent 
to the system of 
estimates $a_1, \ldots, a_n \ll y$, of course with possibly different 
implied constants depending on the field $K$ and the basis. 

We are going to estimate sums of a quadratic character in $\cal O$ by using 
bounds for sums of a quadratic character in $\zb$. The following simple result 
plays a crucial role in this reduction. 
\begin{lem}\label{lem 2.5.2} Let $K$ be any number field and $n$ its degree. 
If the integral ideal $\ff$ of $K$ has square-free norm, then every residue 
class $(\mod \ff )$ is represented by a rational integer.
\end{lem}
\begin{proof}
Let $f = N\ff$. Then, the number of residue classes $\mod \ff$ 
is, by definition, just $f$ and, to show the result, it suffices 
to show that the numbers $1, 2,\ldots, f$ are incongruent
modulo $\ff$. Now, suppose two such rational integers $a$ and $b$ 
are congruent $\mod \ff$. Then, by the complete 
multiplicativity of the norm, 
$N\ff = f$ divides $N(a-b) =(a - b)^n$. But these are rational 
integers and $f$ is square-free so $f$ must divide
the square-free kernel of $(a-b)^n$, which in turn divides $a-b$. 
(Alternatively, ${\cal O}/\,\ff$ is the direct product of prime 
fields $\Fbb_p$ for distinct primes $p$ and so, by the Chinese Remainder 
Theorem, is isomorphic to $\zb/ f\zb$.) 
\end{proof}

Let $\fq$ be an odd ideal. Consider the symbol
\begin{equation}\label{eq:2.15}
\chi_{\fq}(\ell) = \Bigl(\frac{\ell}{\fq}\Bigr) 
\quad {\rm for}\,\, \ell \in \zb\ .
\end{equation}
This is multiplicative in $\ell$ and periodic of period $q = N\fq$ so, as a 
function on $\zb$, $\chi_{\fq}$  is a real Dirichlet character of modulus $q$. 
We need to know when $\chi_{\fq}$ can be the principal character. 
\begin{lem}\label{lem 2.6}
For an odd ideal $\fq$, if $q=N\fq$ is not squarefull then the 
Dirichlet character $\chi_{\fq}$ is not principal. 
\end{lem}
\begin{proof}
Suppose $q$ is not squarefull. That means there is a prime $p$ dividing $q$ 
whose square does not divide $q$. Since $q=N\fq$, there is a prime ideal 
$\fp$ lying above $p$ and dividing $\fq$. Take $\ell$ which is a non-square 
modulo $p$ and $\ell \equiv 1 \, (\mod q/p)$. Such an $\ell$ exists by 
the Chinese Remainder Theorem. Now, $\chi_{\fq}(\ell) 
= \chi_{\fp}(\ell)\chi_{\fq/\fp}(\ell) = -1$ 
so the character is not principal. 
\end{proof}

In the following sections we shall often appeal to some estimates 
for the units in ${\cal U}^+$ which we are now going to present. 
These are not new (see for example, Cassels [Ca], Lang [La], for different 
arguments) but we include brief proofs for completeness.  

\begin{lem}\label{lem 2.2}
There exists a unit $u\in {\cal U}^+$ such that all but one of its conjugates 
are $\le \frac 12$. 
\end{lem}
\begin{proof}
Let $-1,\,\ve_1,\ldots,\ve_r$ be generators of ${\cal U}$ with 
$\ve_{\ell} >0$ for $1\le \ell \le r =n-1$. Take 
$$
u=\ve_1^{a_1}\ldots\ve_r^{a_r}, 
$$
so
$$
u^{(k)}=\ve_1^{(k)a_1}\ldots\ve_r^{(k)a_r}\ .
$$
We need integers $a_1,\ldots, a_r$ such that all the linear forms 
$$
L_k = a_1\log \ve_1^{(k)}+ \ldots + a_r\log \ve_r^{(k)}, \quad 1\le k \le r 
$$
are negative. Because the determinant (the regulator) 
of 
$$
{\cal R} = \bigl( \log \ve_{\ell}^{(k)}\bigr), \quad 1\le \ell,\, k \le r 
$$
does not vanish, there are real numbers $a_1,\ldots, a_r$ with $L_k=-1$ 
for all $1\le k \le r$. Approximate these by rationals and then clear the 
denominators getting a unit $u \in {\cal U}$ with $u^{(1)},\ldots,u^{(r)}<1$. 
Raising $u$ to a sufficiently large even power, we get a unit with
the required property.
\end{proof} 

\begin{lem}\label{lem 2.3}
Let $B\ge 1$. The number of integers $\alpha \in {\cal O}_K$ 
all of whose conjugates satisfy $|\alpha^{(k)}| \le B$ is finite. 
Those of the integers which are units $u \in {\cal U}^+$ also satisfy 
$u^{(k)} \ge B^{-r}$.  
\end{lem}
\begin{proof}
The number of possible irreducible polynomials these integers can satisfy 
is finite since the degree and coefficients of the latter are bounded. 
This gives the upper bound. The lower bound in case of $u \in {\cal U}^+$ 
follows  because 
$Nu=1$.
\end{proof}

Finally, we want to give examples of some fields having the properties we 
have been assuming. The following result is handy for this purpose. 
\begin{lem}
\label{lem 2.9}
Suppose $n$ is an odd prime, and $2$ is a primitive root modulo $n$.  
If $\calu$ contains 
a unit that is neither totally positive nor totally negative, 
then $\calu^+ = \calu^2$.
\end{lem}

\begin{proof}
Let $H = \{\pm1\}^n$, and consider the map 
$\varphi : \calu/\calu^2 \to H$ defined by
$$
u \mapsto (\mathrm{sign}(u^{(1)}),\ldots,\mathrm{sign}(u^{(n)}))\ .
$$
The Galois group $G$ acts naturally on $H$ by 
permuting the real embeddings of $K$.
With this action $H$ is a free $\Fbb_2[G]$-module of rank one, and 
$\varphi$ is an $\Fbb_2[G]$-homomorphism.  If $n$ is an odd prime, 
then $\Fbb_2[G]$ decomposes
$$
\Fbb_2[G] \cong \Fbb_2[x]/(x^n-1) \cong \Fbb_2[x]/(x-1) 
\oplus \Fbb_2[x]/(\Phi_n(x))
$$
where $\Phi_n(x) = (x^n-1)/(x-1)$ is the $n$-th cyclotomic polynomial.  
If $2$ is a primitive 
root modulo $n$, then $\Phi_n(x)$ is irreducible in $\Fbb_2[x]$, and 
we conclude that 
$H$ is the direct sum of two irreducible $\Fbb_2[G]$-modules, the 
one-dimensional trivial representation 
and an irreducible $(n-1)$-dimensional complement.  The image 
$\varphi(\calu)$ 
is an $\Fbb_2[G]$-submodule of $H$ containing 
the trivial subspace $\{\varphi(1),\varphi(-1)\}$, so if $\calu$ 
contains any unit of 
mixed signs then $\varphi(\calu)$ must be all of $H$, i.e., 
$\varphi$ is surjective.

By \eqref{eq:2.2} we have $[\calu:\calu^2] = 2^n = |H|$, so $\varphi$ 
is an isomorphism.  
Since $\calu^+/\calu^2$ is in the kernel of $\varphi$, we conclude 
that $\calu^+ = \calu^2$.
\end{proof}

\noindent{\bf Examples}. 
Suppose $n=3$. 
There is a nice family of cyclic cubic fields, introduced by 
Shanks [Sha], 
which provides examples of number fields satisfying $\calu^+ = \calu^2$. 

For integers $m$, let $\alpha_m$ be a root of the polynomial
$$
f_m(x) = x^3 + mx^2 + (m-3)x - 1\ .
$$
Note that the only rational roots $f_m$ can have are $\pm 1$, 
but $f_m(1) = 2m-3$ and $f_m(-1) = 1$, so $f_m$ is irreducible for every 
$m$. The discriminant of $f_m$ is $(m^2 - 3m + 9)^2$, so 
$K_m = \qb(\alpha_m)$ is a cyclic cubic field for every $m$.  

Next, we note that for every $m$, Descartes' ``Rule of Signs'' shows that 
the polynomial $f_m(x)$ has exactly one positive real root.  
Therefore $\alpha_m$ 
has both positive and negative real embeddings, so Lemma \ref{lem 2.9} 
shows that $\calu^+ = \calu^2$ for every integer $m$.

Similarly, when $n=5$ there is a family of cyclic 
quintic fields constructed by 
Lehmer [Le].  For every integer $m$ let 
\begin{equation*}
\begin{split}
g_m(&x) = x^5 + m^2x^4  - 2(m^3+3m^2+5m+5)x^3 \\
 &+(m^4+5m^3 +11m^2 +15m +5)x^2 
 +(m^3+4m^2+10m+10)x+1\ ,
\end{split}
\end{equation*}
and let $\beta_m$ be a root of $g_m$.  The field 
$L_m = \qb(\beta_m)$ is a cyclic quintic 
extension of $\qb$, and $\beta_m$ is a unit of $L_m$ of norm $-1$.  
Thus $\beta_m$ has at least one 
negative real embedding.  If all the real embeddings of $\beta_m$ 
were negative, then 
all coefficients of $g_m$ would be positive.  It is simple to 
check that 
$-2(m^3+3m^2+5m+5)$ and $m^3+4m^2+10m+10$ are not simultaneously 
positive for any integer $m$, so 
$\beta_m$ has at least one positive real embedding.
Since $2$ is a primitive root modulo $5$, we conclude from 
Lemma \ref{lem 2.9} that 
$\calu^+ = \calu^2$ for every integer $m$.

\section{A Fundamental Domain}\label{sec 3}

In order to have a convenient unique representation of a principal ideal 
by one of its generators, we look for a specific fundamental domain of the 
group $\cal U^+$ acting on $\rb_+^n$ by
\begin{equation}\label{eq:3.1}
u\circ x = \bigl(u^{(1)}x_1,\ldots, u^{(n)}x_n\bigr)\ . 
\end{equation}
For notational simplicity, we write $x\succ 0$, meaning that 
all the co-ordinates 
of $x=(x_1,\ldots,x_n)$ are positive. Simlilarly, $x\succ C$ or $x\prec C$ 
 means that 
all the co-ordinates are greater than $C$, or smaller than $C$, respectively. 
For $x=(x_1,\ldots,x_n)$ and $y=(y_1,\ldots,y_n)$ in  $\rb^n$ the scalar 
product is defined by
$$
x\cdot y = x_1y_1 +\ldots + x_ny_n\ .
$$

We begin with a somewhat general consideration. Let $U$ be a collection of 
$n$-dimensional positive vectors which does not contain $e=(1,\ldots,1)$. 
Put
\begin{equation}\label{eq:3.2}
{\cal D} =\{x\in \rb_+^{n}\,; \,\,\,x\succ 0\, ,\quad u\cdot x>e\cdot x
\quad {\rm for\,\, all}\,\, u\in U\}\ . 
\end{equation}

Suppose there is a matrix $(c_{k\ell}), \, 1\le k,\,\ell \le n$ whose rows  
are in $U$, with its diagonal entries $c_{kk} > 1$ for $1\le k\le n$ 
and its off-diagonal entries satisfying 
$0<c_{k\ell}<1$ for $1\le k \neq \ell\le n$. Put 
\begin{equation}\label{eq:3.3}
C= 1 + \max_{k\neq \ell} \frac{c_{kk}-1}{1-c_{k\ell}}\ .
\end{equation}
Note that all the entries of $(c_{k\ell})$ satisfy 
$0< c_{k\ell} < C$. Define 
\begin{equation}\label{eq:3.5}
{\cal D} (C) 
= \{x\in \rb_+^{n}\,; \,\,\,x\succ 0\, ,\,\, v\cdot x> e\cdot x\,\,\, 
{\rm for\,\, all} \,\, v \in U,\, v\prec C\}\ .
\end{equation}
\begin{lem}\label{lem 3.1}
We have
\begin{equation}\label{eq:3.6}
{\cal D}(C) = \cal D \ .
\end{equation}
\end{lem}
\begin{proof}
Obviously, ${\cal D} \subset {\cal D} (C)$. Suppose there exists 
$x\in {\cal D} (C)\setminus \cal D$. Because $x\not\in {\cal D}$ there 
exists $u\in U$ such that $u\cdot x\le e\cdot x$. 
Hence, $u$ is not bounded by $C$ 
because $x\in {\cal D} (C)$. This means that $u=(u_1,\ldots,u_n)$ has 
$u_k\ge C$ 
for some $ 1\le k\le n$. We write the inequality $u\cdot x\le e\cdot x$ for 
$x=(x_1,\ldots,x_n)$ in the following form:
$$
(u_k-1)x_k \le \sum_{\ell \ne k}(1-u_{\ell})x_{\ell}\ .
$$
Hence,
$$
(C-1)x_k < \sum_{\ell \ne k}x_{\ell}\ .
$$
But 
$$
(C-1)\ge (c_{kk}-1)/ \min_{\ell\neq k}(1-c_{k\ell})
$$
and hence 
$$
(c_{kk}-1)x_k <\sum_{\ell \ne k}(1-c_{k\ell})x_{\ell}\ .
$$
Since the vector $v=(c_{k1},\ldots,c_{kn})$ belongs to $U$ 
and its entries are bounded by $C$, the last inequality contradicts 
the assumption that $x\in {\cal D} (C)$. 
\end{proof}

We are going to apply Lemma~\ref{lem 3.1} to the collection 
of vectors 
$$
U= \bigl\{(u^{(1)},\ldots,u^{(n)});\,\, u\in {\cal U}^+,\, u\neq 1\bigr\}\ .
$$
By Lemma~\ref{lem 2.2} there exists a unit $u\in {\cal U}^+$ such that 
$u^{(1)},\ldots,u^{(n-1)}\, <1$ and hence $u^{(n)} >1$. This unit $u$ and its 
conjugates form the matrix $(c_{k\ell})$ whose existence was used 
in the proof of Lemma~\ref{lem 3.1}. We obtain:
\begin{lem}\label{lem 3.2} 
There exists a constant $C>1$ such that the sets
\begin{equation}\label{eq:3.7}
{\cal D}  
= \{ x\succ 0;\,\, u\cdot x> e\cdot x\,\,\, {
\rm for\,\, all} \,\, u \in {\cal U}^+,\, u\neq 1\}\ ,
\end{equation}
\begin{equation}\label{eq:3.8}
{\cal D} (C) 
= \{ x\succ 0;\,\, u\cdot x> e\cdot x\,\,\, 
{\rm for\,\, all} \,\, u \in {\cal U}^+,\, u\neq 1, \, u \prec C\}\ , 
\end{equation}
are the same set.
\end{lem}
Note that, by Lemma~\ref{lem 2.3}, the collection of units 
$u\neq 1,\,u\prec C$ 
which are used in ${\cal D}(C)$ is finite. We fix such a $C$ and 
denote 
\begin{equation}\label{eq:3.7.1}
\tilde{\cal U} = \{ u\in {\cal U}^+; \, u \neq 1, \, u\prec C\}\ .  
\end{equation}
Note that the closure ${\bar{\cal D}}$ consists of those $x$ for which 
all of the signs $>$ in $\quad u\cdot x>e\cdot x$ 
are replaced by $\ge$ and the boundary  
consists of those elements of the closure for which we have equality 
with at least one unit $u \neq 1$. 

Now, we are ready to prove the fundamental properties of the 
set~\eqref{eq:3.7}. 
\begin{lem}\label{lem 3.3} 
For any $u\in {\cal U}^+,\, u\neq 1$ we have 
\begin{equation}\label{eq:3.8.1}
(u\circ {\cal D})\cap {\cal D} = \varnothing\ . 
\end{equation}
For every $x>0$ there exists $u\in {\cal U}^+$ which sends $x$ to 
the closure of ${\cal D}$, that is
\begin{equation}\label{eq:3.8.2}
u\circ x \in \bar{\cal D}\ . 
\end{equation}
\end{lem}

\begin{proof}
Suppose $x\in {\cal D}$ and $x\in u\circ{\cal D}$. Then,
$y=u^{-1}\circ x \in {\cal D}$, so $u^{-1}\cdot x>e\cdot x$ and $u\cdot y
>e\cdot y$. The second 
inequality reads as 
$$e\cdot x= u\cdot (u^{-1}\circ\,x) > e\cdot (u^{-1}\circ\,x)  
= u^{-1}\cdot x\ ,
$$
 which contradicts the first 
inequality. This proves the first property~\eqref{eq:3.8.1}. 
The second property~\eqref{eq:3.8.2} follows by choosing $u\in {\cal U}^+$ 
for which $u\cdot x$ is minimal. The minimum is 
attained since it is over a finite set, which is seen to be the case 
because ${\cal D} = {\cal D}(C)$. Suppose this $u$ did not have the 
property~\eqref{eq:3.8.2}. This means that there is a $v\in {\cal U}^+$, 
$v\neq\, {\rm id}$, such that $v\cdot (u\circ x)<u\cdot x$. Hence 
$vu \in \cal U^+$ gives 
a smaller value than the minimum and proves~\eqref{eq:3.8.2}.
\end{proof} 

The above argument is essentially our interpretation of an argument of 
Shintani [Shi] who gave a complete description of the fundamental domain, say 
${\cal D}^*$, of the action of ${\cal U}^+$ on $\rb_+^n$. Lemma~\ref{lem 3.3} 
shows that one can choose ${\cal D}^*$ with
\begin{equation}\label{eq:3.9}
{\cal D}\subset {\cal D}^* \subset \bar{\cal D}\ . 
\end{equation}

For our applications we do not need to see the boundary of ${\cal D}^*$ 
and that simplifies things a lot. 
Note that, if $\alpha\in{\cal O}$ then $\alpha\in{\cal D}$ 
means that $\alpha \succ 0$ and it has the smallest trace amongst 
its positive associates $u\alpha, \, u\in{\cal U}^+$. 
There are several features of this particular domain ${\cal D}^*$, 
such as convexity, which help to control the transition from 
ideals to integers. We present a few of these here; the others will be 
introduced as they are exploited in Section~\ref{sec 5}. 

\begin{lem}\label{lem 3.4}
Every $\alpha \in {\cal D}$ has all its conjugates in ${\cal D}$ 
and all of them have essentially the same 
size, that is 
\begin{equation}\label{eq:3.10}
\alpha^{(k)}\asymp T(\alpha)\asymp (N\alpha)^{\frac 1n}, \quad 1\le k\le n\ . 
\end{equation}
\end{lem} 
\begin{proof}
Choose $u\in {\cal U}^+$ such that $u^{(\ell)} \le \frac 12$ for all 
$\ell \neq k$, as in Lemma~\ref{lem 2.2}. Then, we see in turn,
\begin{equation*}
\begin{split}
u^{(1)}\alpha^{(1)}+\ldots +u^{(n)}\alpha^{(n)} & >
\alpha^{(1)}+\ldots + \alpha^{(n)} , \\  
u^{(k)}\alpha^{(k)} & > \tfrac 12 \bigl(\alpha^{(1)}+\ldots + \alpha^{(n)}
\bigr) , \\  
\alpha^{(k)} & \gg T\alpha, \quad k=1,2,\ldots ,n\ .
\end{split}
\end{equation*}
Now, take 
$u\in {\cal U}^+$ to be the reciprocal of the previous choice, giving 
$\alpha^{(k)}\ll T\alpha$.  
This completes the proof. 
\end{proof}
\begin{cor}\label{cor 3.5}
If $\alpha = a_1\omega_1 + \ldots +a_n\omega_n \in {\cal D}$ then 
\begin{equation}\label{eq:3.11}
a_k \ll (N\alpha)^{\frac 1n}\quad {\rm for}\,\, k =1, \ldots, n\ .
\end{equation}
\end{cor}

Denote by ${\cal N}(x)$ the number of integers $\alpha \in \bar{\cal D}$ 
with $N\alpha\le x$. Thus 
\begin{equation}\label{eq:3.12}
{\cal N}(x)\asymp x\ ,
\end{equation}
but we do not need an asymptotic formula for ${\cal N}(x)$. More important 
to us is to have a precise comparison of ${\cal N}(x)$ 
with the number of these integers which are in a given residue class.  
\begin{equation}\label{eq:3.13}
{\cal N}(x;\fm,\nu)= \sum_{\substack{\alpha\in \bar{\cal D},\,N\alpha\le x\\ 
\alpha \equiv \nu\,(\,\mod \fm)}}1\ . 
\end{equation}
\begin{lem}\label{lem 3.6} For any integral ideal $\fm$ and any 
$ \nu\,(\mod \fm) $ we have
\begin{equation}\label{eq:3.14}
{\cal N}(x;\fm,\nu)= \frac{{\cal N}(x)}{N\fm} +O\bigl(x^{1-\frac 1n}\bigr)\ ,
\end{equation}
where the implied constant depends only on the field $K$ (and so, 
not on $\fm$). 
\end{lem} 

\begin{proof} Fix an integral basis of $\cal O$ and choose one of its 
elements, say $\omega$. 
Write $\alpha = a\omega +\beta$ where $a\in \zb$ and $\beta$ is a linear  
combination of the other basis elements. 
Given $\beta$,  
we are going to evaluate the number of rational integers $a$ such that 
\begin{equation}\label{eq:3.15}
a\omega + \beta \in \bar {\cal D},\quad P(a)\le x\ , 
\end{equation}
where $P(X) = \bigl(X\omega^{(1)}+\beta^{(1)}\bigr)
\ldots\bigl(X\omega^{(n)}+\beta^{(n)}\bigr)$, and 
\begin{equation}\label{eq:3.16}
a\omega + \beta \equiv \nu\, (\mod \fm)\ . 
\end{equation}
The first condition in~\eqref{eq:3.15} means that $a$ satisfies 
a finite system of linear inequalities which also means exactly that $a$ 
is in a single interval whose end-points depend on $\omega, \, \beta$ and $K$. 
The second condition in~\eqref{eq:3.15} can be expressed as saying 
that $a$ runs over $n$ intervals, each of length $\ll x^{1/n}$, 
whose end-points depend on $\omega, \, \beta$ and $x$. Together, 
the two conditions in~\eqref{eq:3.15} are equivalent to saying that $a$ 
is in one of $n$ intervals of length $\ll x^{1/n}$ 
whose end-points depend on $\omega, \, \beta, \, K$ and $x$.

Next, the congruence condition~\eqref{eq:3.16}, if solvable, means that
we have $a\equiv a_0 \, ( \mod m)$ where $m$ is a positive rational 
number which 
depends on $\omega, \, \fm$ and $a_0$ is a rational residue class modulo 
$m$ which depends on $\omega, \, \nu, \, \beta$ and $\fm$. 
(We can take $m$ to be the smallest 
positive rational integer divisible by the ideal $\fm/(\omega, \fm)$.) 
Changing $\nu$ to $\nu' = \nu + c\,\omega$ with $c\in \zb$ 
translates the class $a_0$ to $a_0 +c$. 
This operation can change the number of $a$'s in a given segment 
(of an arithmetic progression) by at 
most one. Such a bounded error term is then amplified by the number of 
segments and by the number of $\beta$ with $N\beta \ll x$ which is 
$O(x^{1-1/n})$. This proves that 
$$
{\cal N}(x;\fm,\nu)={\cal N}(x;\fm,\nu') +O\bigl(x^{1-\frac 1n}\bigr)\ ,
\quad {\rm if} \,\, \nu' \in \nu + \omega\zb\ ,
$$
where the implied constant depends only on the 
field $K$.
Repeating these arguments with every element of the basis 
we derive the result for every $\nu' (\mod \fm)$.  
Then, averaging this relation over all classes $\nu' (\mod \fm)$,
we complete the proof of~\eqref{eq:3.14}. 
\end{proof}

\section{Sums over Prime Ideals}\label{sec 4}

We begin with a few formulas of a combinatorial nature. 
Let $K$ be a number field. For any non-zero integral ideal $\fn$ we set
$$
\Lambda (\fn)= \log N\frak p\ ,\,\, {\rm if}\, \fn =\frak p^{\ell}, 
\quad \ell = 1,2,\ldots\ ,
$$
and $\Lambda (\fn)= 0$ otherwise.
Hence it is easy to see that
$$
\sum_{\fb \mid \fa} \Lambda (\fb) = \log N\fa\ . 
$$
Next, we introduce the M\"obius function 
$$\mu (\fm) = (-1)^t
$$
if $\fm$ is the product of $t$ distinct prime ideals and $\mu (\fm)=0$ 
otherwise. 
Note that for $\fm =(1)$ we have $t=0$ and $\mu ((1))=1$. Hence, we deduce    
$$
\sum_{\fm \mid \fa} \mu(\fm)  =  \, 
\begin{cases}\, 1 
 \quad \text{if $\fa = (1)$}\, , 
\\
0 
 \quad \text{otherwise}\, .
\end{cases}
$$
Using this, one can check the formulas
$$
\Lambda (\fn) = - \sum_{\fm \mid \fn} \mu(\fm)\log N\fm
= \sum_{\fm \mid \fn} \mu(\fm)\log \frac{N\fn}{N\fm}\ . 
$$

Let ${\cal A}= ( a_{\frak n})$ be an arbitrary sequence of complex numbers,  
enumerated by integral ideals and ordered by the norm. 
We are interested in estimating the sum
\begin{equation}\label{eq:4.1} 
S(x) = \sum_{N\frak n \le x}a_{\frak n}\Lambda(\frak n)\ . 
\end{equation} 
We are thinking of the $a_{\frak n}$ as changing argument randomly and 
expect considerable cancellation in the sum $S(x)$. Having this in mind, 
we are going to write a fairly general inequality which offers a bound 
for $S(x)$ in terms of other sums which we know how to manage. The idea 
goes by adding more terms in the spirit of the Eratosthenes-Legendre sieve 
until reaching two kinds of sums. The first kind are the congruence sums 
\begin{equation}\label{eq:4.2} 
A_{\frak d}(x) = \sum_{\substack{N\frak n \le x\\ 
\frak n\equiv 0\, (\mod\, \frak d)}}a_{\frak n}\ . 
\end{equation} 
These will appear for $\frak d$ with $d= N\frak d$ relatively small, so 
the problem of estimating $A_{\frak d}(x)$ really belongs to the harmonic 
analysis of ${\cal A} = (a_{\frak n})$. 

The second kind of sums are the bilinear forms 
\begin{equation}\label{eq:4.3} 
{\cal B} (M,N) = \sum_{N\frak m \le M}\sum_{N\frak n \le N}
v_{\frak m}\, w_{\frak n}\, a_{\frak m \frak n}\ . 
\end{equation} 
These will appear for $M, \, N$ neither of which is very small so that
${\cal B} (M,N)$ is a genuine bilinear form. Here the point is that 
$(v_{\frak m}), \, (w_{\frak n})$ are independent sequences; 
they do not see each other so they cannot conspire to annihilate the 
change of arguments of $a_{\frak m \frak n}$. Well, except for a sequence 
$\cal A$ whose terms are multiplicative; for example, if 
$a_{\frak m \frak n}= c_{\frak m}c_{\frak n}$ then the bilinear 
form~\eqref{eq:4.3} factors into linear forms 
$$
{\cal B} (M,N) = \Bigl(\sum_{N\frak m \le M}
v_{\frak m}\,c_{\frak m}\Bigr)\,\,
\Bigl(\sum_{N\frak n \le N} w_{\frak n}\, c_{\frak n}\Bigr)  
$$
and we can obtain a bias by choosing $v_{\frak m}= \bar c_{\frak m}$
and $w_{\frak n} = \bar c_{\frak n}$. Therefore, it is important that the 
sequence ${\cal A} = (a_{\frak n})$ not be multiplicative in $\frak n$. 
Our target sequence $a_{\frak n} = spin(\frak n)$ is qualified for 
treatment by our method, due 
to the twisted factorization property~\eqref{eq:2.10}. 

We shall see that the required bilinear forms ${\cal B} (M,N)$ have specific 
coefficients $v_{\frak m}, \, w_{\frak n}$, but in practice we 
are unable to take advantage of their intrinsic properties, so there is 
no point to describe these. Our estimates for ${\cal B} (M,N)$ will depend 
only on the upper bound for their coefficients. Specifically, we shall 
use the bilinear form with 
\begin{equation}\label{eq:4.4} 
|v_{\frak m}| \le \Lambda (\frak m), \quad\quad 
|w_{\frak n}| \le \tau(\frak n)\ ,
\end{equation} 
where $\tau$ is the usual divisor function, but of ideals 
in the ring ${\cal O}$.

\begin{prop}\label{prop 4.1}
Let $x=yz$ with $z\ge y\ge 2$. We have 
\begin{equation}\label{eq:4.5} 
|S(x)| \le (3\log x)\sum_{N\frak d\le y^2}|A_{\frak d}(x')| 
+ \,(2\log x)^2 \, |{\cal B}(M,N)| + |S(z)|\ , 
\end{equation} 
for some $x'\le x$, some $M,\, N \le z$, $MN=2x$ and some complex coefficients 
$v_{\frak m}\, , \, w_{\frak n}$ in ${\cal B} (M,N)$ 
satisfying~\eqref{eq:4.4}.
\end{prop}

\begin{proof} 
We begin by decomposing the convolution $\Lambda = \mu * \log$ as follows:
\begin{equation*}
\Lambda (\fn) 
= \sum_{\fa\fm =\fn}\mu(\fm)\log N\fa 
= \sum_{\substack{\fa\fm =\fn\\ N\fm \le y}}\mu(\fm)\log N\fa \,\, 
+ \, \sum_{\substack{\frak l \fa\fm =\fn\\ N\fm > y}}\mu(\fm)\Lambda (\fa)\ . 
\end{equation*}
Suppose $\fn$ has norm $N\fn \le x =yz$ with $z\ge y \ge 2$. Then, in 
the last sum we have $N\fa \le z$. Having recorded this information, we 
now write
$$
\sum_{\substack{\frak l \fa\fm =\fn\\ N\fa\le z\, ,N\fm > y}}
\mu(\fm)\Lambda (\fa) =
\sum_{\substack{\frak l \fa\fm =\fn\\ N\fa\le z}}\mu(\fm)\Lambda (\fa)\,\,\, 
- \sum_{\substack{\frak l \fa\fm =\fn\\ N\fa\le z\, ,N\fm  \le y}}
\mu(\fm)\Lambda (\fa)\ .
$$
If, in the first sum on the right, we fix $\frak a$, the 
inner {\it complete} sum over $\frak l \fm$ vanishes 
unless $\frak l \fm =(1)$ in which case $\fa =\fn$. Hence, we get the 
following identity:
$$
\bigl(1 -\delta(\frak n ,z)\bigr)\Lambda (\frak n) = 
\ssum_{\substack{\frak a \frak m=\frak n\\N\frak m \le y}}
\mu (\fm)\log N\frak a \,\,
-\,\, \sssum_{\substack{\frak l \frak a \frak m = \frak n\\N\frak a\le z, \,
N\frak m \le y}}\mu (\frak m) \Lambda (\frak a)\ ,  
$$
where $\delta(\frak n, z)=1$ if $N\frak n \le z$ and is zero 
elsewhere. We split the sum over $\frak l \frak a\frak m =\frak n$ into 
two sums having $N\frak a \le y$ or not. Accordingly, $S(x)-S(z)$ splits 
into three sums 
$$
S(x)-S(z) = S_1(x)-S_2(x)-S_3(x)\ ,
$$ 
where 
$$
S_1(x) = \sum_{N\frak m\le y}\mu(\frak m)\sum_{\substack{N\frak n \le x\\ 
\frak n \equiv 0\,(\mod \frak m )}}
a_{\fn} \log N\bigl(\frac{\frak n}{ \frak m} \bigr)\ ,
$$
$$
S_2(x)= \sum_{\frak d}\Bigl(\sum_{\substack{\frak a \frak m =\frak d\\ 
N\frak a\le y,\, N\frak m\le y}}\mu(\frak m)\Lambda (\frak a)\Bigr)
A_{\frak d}(x) 
$$
and
$$
S_3(x)= \sssum_{\substack{N (\frak l \frak a \frak m )\,
 \le x \\ y\, < N\frak a\, \le z, \,
N\frak m \, \le y}}\mu (\frak m) \Lambda (\frak a)
a_{\frak l \frak a \frak m}\ . 
$$ 
Note that, in the sum $S_2(x)$, we have $N\frak d \le y^2$ and 
the coefficient in front of the congruence sum 
$A_{\frak d}(x)$ is bounded by
$$
\sum_{\frak a \frak m = \frak d}\Lambda (\frak a) = \log N\frak d \ . 
$$
Similarly, we treat the first sum $S_1(x)$.  
Here, the presence of $\log N (\frak n / \frak m )$ is somewhat inconvenient 
so we replace it by 
$$
\log N\bigl(\frac{\frak n}{ \frak m}\bigr) 
= \int^{N\frak n}_{N\frak m}t^{-1}dt
$$
 and then, inverting the order of 
summation and integration, we arrange $S_1(x)$ into
\begin{equation*}
\begin{split}
S_1(x) & = \sum_{N\frak m\le y}\mu(\frak m)
\int^{x}_{N\fm}\bigl(A_{\frak m}(x)-A_{\frak m}(t\bigr))t^{-1}dt\\
& = \int^{x}_{1}\sum_{N\frak m\le \min (y, t)}\mu(\frak m)
\bigl(A_{\frak m}(x)-A_{\frak m}(t\bigr))t^{-1}dt\ . 
\end{split}
\end{equation*} 
Now, we bound the inner sum by
$$
\max_t \big|\sum_{N\frak m\le \min (y, t)}\mu(\frak m)
\bigl(A_{\frak m}(x)-A_{\frak m}(t\bigr))\big|
$$
and pull it outside the integral. We deduce that
$$
|S_1(x)| \le 2\sum_{N\frak m\le y}|A_{\frak m}(x')|
\log x\ , 
$$
for some $x'\le x$. 
Adding the bound for $S_2(x)$ we get 
$$
|S_1(x)|  + |S_2(x)| \le 3(\log x)\sum_{N\frak d\le y^2}|A_{\frak d}(x')|
$$
for some $x'\le x$. 

We consider the triple sum $S_3(x)$ 
as a double sum with $\frak a$ being 
one variable and $\frak b = \frak l \frak m$ as the second variable. These 
variables are weighted by $v_{\frak a}=\Lambda (\frak a)$ and 
$$
w_{\frak b} = \sum_{\frak m | \frak b, \, N\frak m \le y}\mu(\frak m), 
\quad {\rm so}\, |w_{\frak b}| \le \tau(\frak b)\ . 
$$
The variables are restricted by $y < N\frak a \le z$ and 
$N\frak a \frak b \le x$. Hence, $N\frak b \le xy^{-1} =z$. Now we 
would like to relax the condition $N\frak a \frak b \le x$ because it ties 
the two variables together, but still we should be able to recover a slightly 
weaker condition $N\frak a \frak b \le 2x$. To this end, we subdivide 
the range for $N \frak a$ into dyadic segments $\frac 12 M< N\frak a \le M$ 
starting with $M=z$. The number of such intervals is no more than 
$\log z /\log 2$. 
When $a=N\frak a$ is in such an interval then $b=N\frak b \le x/a<2x/M=N$, 
say. Having recorded that $a\le M, \, b\le N$ with $M,\, N\le z, \, MN=2x$, 
we remove the condition $ab\le x$ by a standard technique of separation 
of variables.
Lemma 9 of [DFI] provides us with a function $h(t)$ which satisfies 
$\int^{\infty}_{-\infty} |h(t)|\, dt \le \log 6x$ and also, 
for positive integers $k$, 
$$
\int^{\infty}_{-\infty} h(t) k^{it} dt  =  \, 
\begin{cases}\, 1 
 \quad \text{if $1\le k \le x$}\, , 
\\
0 
 \quad \text{if $k > x$}\, .
\end{cases}
$$
We insert this integral, with $k=ab$, as a factor in the summation,  
allowing us to, in effect, separate the variables $a$ and $b$. 
We then interchange summation and integration. In this inner summation 
the coefficients of $\frak a$ and $\frak b$ are now 
contaminated by the twists $a^{it}$ and $b^{it}$ for a real $t$. 
 These 
contaminating factors have absolute value one, hence do not change the bounds 
for the coefficients $v_{\frak a}, \, w_{\frak b}$. The
summation now depends on $t$, but we majorize by choosing 
that $t$ which maximizes the absolute value of the whole sum. 
Having rendered it independent of $t$, we now pull this absolute value 
outside the integral, 
then integrate. The integration costs us a factor 
$\log 6x$ and, in total, the above 
operations cost us a loss in the bounds for the coefficients by a factor 
$(\log z /\log 2) \log 6x \le (\log x /\log 2) \log 6x \le (2\log x)^2$. 
This completes the proof of~\eqref{eq:4.5}. 
\end{proof} 

One can show that~\eqref{eq:4.5} holds with the congruence sums 
$A_{\frak d}(x)$ restricted by $N\frak d \le y$. This little improvement  
can be achieved along the above lines by more careful partitions, however, 
such a result would have no significance in our applications. 

\begin{prop}\label{cor 4.2}
Suppose we have fixed numbers $0<\vartheta,\,\theta<1$ such that 
the sequence ${\cal A} = (a_{\frak n})$ with $|a_{\frak n}| \le 1
$ allows the following estimations: 
\begin{equation}\label{eq:4.7} 
A_{\frak d}(x) \ll  x^{1-\vartheta +\ve}\ , 
\end{equation} 
for any ideal $\frak d$ and  
any $x\ge 2$, and
\begin{equation}\label{eq:4.8} 
{\cal B}(M,N) \ll (M+N)^{\theta}(MN)^{1-\theta +\ve}\ , 
\end{equation}
for any $M,\, N \ge 2$. Here $\ve$ is any positive number 
and the implied constants depend
only on $\ve$ and the field $K$. Then, for any $x\ge 2$, we have 
\begin{equation}\label{eq:4.9} 
S(x) \ll x^{1- \frac{\vartheta\theta}{2+\theta} +\ve} 
\end{equation}
with any $\ve >0$, the implied constant depending on $\ve$ and the field $K$. 
\end{prop} 
\begin{proof}
Use the bounds $S(z)\ll z=xy^{-1}, \,\, |\{ \fd ; N\frak d \le y^2\}| \ll y^2$
and apply~\eqref{eq:4.5} for $y= x^{\vartheta/(2+\theta)}$. 
\end{proof}

We shall verify that~\eqref{eq:4.7} holds for our sequence of spins in 
Section~\ref{sec 5} and~\eqref{eq:4.8} in 
Section~\ref{sec 6}. 

\section{Congruence Sums in the Spin}\label{sec 5}

Recall that $\fM$ is the modulus occurring in the statement of 
Theorem~\ref{theo 1.2}. Since we are targeting ideals in a progression 
modulo $\fM$ it is convenient to introduce the characteristic function 
on these. We define the function $r(\fa) = r(\fa ; \fM, \mu)$, on all 
integral ideals $\fa$ of $K$ by setting 
$r(\fa) =1$ if there exists an integer $\alpha$ of $K$ such 
that $\fa = (\alpha), \, \alpha \succ 0, \, \alpha \equiv \mu\, (\mod\, \fM)$  
and we put $r(\fa) =0$ otherwise.
Keep in mind that $r(\fa)$ is supported on ideals 
co-prime with $\fM$, so on odd ideals. Our final goal is to sum  
the spin of those primes in the support of $r$.  

Let $F$ be a fixed positive rational integer which is a 
multiple of $\fM$ and of $2^{2h+3}f$, where 
$f$ is given by~\eqref{eq:2.9} and $h$ is the class number.  
Let $\frak m$ be an ideal, co-prime with its Galois conjugate 
$\frak m'$ and with $F$. We consider 
\begin{equation}\label{eq:5.1} 
A(x)= \sum_{\substack{N\frak a\, \le x\\
(\frak a, F\,)=1,\,\, \frak m | \frak a}} r(\fa)\, spin\, (\frak a)\ .
\end{equation}
We suppress, in the notation, the dependence of the congruence 
sum $A(x)$ on $\fM, \, \mu,\, F$ and 
$\fm$; however, in the estimations we shall pay attention to uniformity 
in terms of $\fm$, but not on the other three which are fixed for us. 
Obviously, $A(x)$ is bounded by the number of all ideals divisible by 
$\frak m$ and having norm $\le x$ so that
\begin{equation}\label{eq:5.2} 
A(x) \ll \frac{x}{N \frak m}\ . 
\end{equation}
Our goal in this section is to prove a much stronger bound for small $N\fm$,  
by exploiting the cancellation due to the sign change of the spin. 
\begin{prop}\label{prop 5.1}
Assume Conjecture $C_n$ (see\eqref{eq:8.4}). Then, for any $\fm$ 
with $(\fm, \fm' F)=1$ and
any $x\ge 2$ we have
\begin{equation}\label{eq:5.3} 
A(x) \ll x^{1-\frac{\delta}{2n} +\ve}\ , 
\end{equation} 
where the implied constant depends on $\ve$ and the field $K$. 
The bound~\eqref{eq:5.3} 
with exponent $\delta = 1/48$ holds unconditionally 
for cubic fields. 
\end{prop}
We remark that the latter statement follows from the former because, for $n=3$ 
the conjecture $C_3$ holds true by Burgess' theorem, with exponent 
$\delta = 1/48$; see Corollary~\ref{cor 8.1}. 

We begin the proof by picking the unique generator of $\fa$, say 
$\fa = (\alpha)$, with $\alpha\in {\cal D}^*$. Recall that, according to 
our convention, $\alpha\in {\cal D}^*$ means
\begin{equation}\label{eq:5.6} 
\bigl(\alpha^{(1)},\ldots,\alpha^{(n)} \bigr)\in {\cal D}^*\ . 
\end{equation} 
Here ${\cal D}^*$ denotes the fundamental domain of the group of 
totally positive units ${\cal U}^+={\cal U}^2$ acting on $\rb_+^n$ as in 
Section~\ref{sec 3}. We do not need to know exactly what the boundary  
$ {\cal D}^*\setminus{\cal D}$ looks like because the contribution of
$\alpha\in {\cal D}^*\setminus{\cal D}$ is negligible for our purpose. 
This unique generator $\alpha$ may not be the same as the one  
equivalent to
$\mu$ modulo $\fM$ whose existence is implied by the support of 
$r(\fa)$, so we claim only that $\alpha \equiv \mu\, u$ for some 
$u \in {\cal U}^+$. Obviously, $u$ is determined up to the units in 
$$
{\cal U}^+_{\fM} = \{v \in {\cal U}^+; \, v \equiv 1\, (\mod \fM)\}\ . 
$$ 

We split the sum~\eqref{eq:5.1} over $\alpha\in {\cal D}$ into residue 
classes modulo $F$ getting
\begin{equation}\label{eq:5.7} 
A(x)=\sumwedge_{\substack{\rho\,(\mod F)\\
(\rho,\, F)=1}} A(x; \rho) +\partial A(x)\ ,
\end{equation} 
where the superscript $\wedge$ restricts the summation to classes 
$\rho \, (\mod F), \,  \rho \equiv \mu\, u\, (\mod \fM)$ 
for some $u \in {\cal U}^+$,  
\begin{equation}\label{eq:5.8} 
A(x; \rho)= \sum_{\substack{\alpha\in {\cal D},\, N\alpha\le x\\
\alpha\equiv \rho\,(\mod F)\\ \alpha\equiv 0 \,(\mod \fm)}}spin (\fa)\ , 
\end{equation} 
and $\partial A(x)$ denotes the contribution of the boundary terms, so 
$$
|\partial A(x)| \le |\{\alpha\in \bar{\cal D}\setminus{\cal D}\,;\, 
N\alpha \le x\}|\ .
$$

It is easy to estimate $\partial A(x)$ so we do it now. 
The condition $\alpha\in \bar{\cal D}\setminus{\cal D}$ implies that there 
is a unit $u\neq 1$ in the finite set $\tilde{\cal U}\subset {\cal U}^+$  
defined in~\eqref{eq:3.7.1},  
such that 
$$
u^{(1)}\alpha^{(1)} +\ldots + u^{(n)}\alpha^{(n)} 
= \alpha^{(1)} +\ldots + \alpha^{(n)}\ ;
$$
see~\eqref{eq:3.7.1} and the remarks following it. 

Returning to the 
notation of Section~\ref{sec 2}, we consider the fixed integral basis 
$1, \omega_2,\ldots,\omega_n$ of $\cal O \subset K$. 
Writing $\alpha = a_1 +a_2\omega_2+ \ldots +a_n\omega_n$ in terms of that 
basis,
the above equation becomes 
$$
\bigl(u^{(1)}+ \ldots +u^{(n)}- n\bigr)\, a_1 
= \lambda_2a_2 +\ldots + \lambda_na_n 
$$
for certain explicit $\lambda_i \in {\cal O}$, independent of $\alpha$. 
Since the factor 
in front of $a_1$ is positive (see~\eqref{eq:2.3}) for $u\neq 1$, the 
coordinate $a_1$ is determined by the other coordinates $a_2,\ldots,a_n$. 
However, all the basis co-ordinates are 
$\ll x^{1/n}$ (see~\eqref{eq:3.11}), so
\begin{equation}\label{eq:5.9} 
\partial A(x)\ll x^{1-\frac 1n}\ ,
\end{equation} 
where the implied constant depends only on the field $K$. Note that we 
have abandoned the condition $\alpha\equiv 0\,(\mod \fm)$ losing something 
in the bound~\eqref{eq:5.9}. This is not a serious issue because the 
estimates to come will be much weaker than~\eqref{eq:5.9}, anyway. 

Next, we are going to estimate every $A(x; \rho)$ separately for the 
residue classes $\rho \, (\mod F),\,(\rho,F)=1$. Returning again to the 
notation of Section~\ref{sec 2}, we consider the module 
$$
\mb = \omega_2\,\zb + \ldots + \omega_n\,\zb 
$$
with rank $\mb = n-1$ and ${\cal O} = \zb +\mb$. 
We now write $\alpha$ uniquely as
\begin{equation}\label{eq:5.10} 
\alpha = a+ \beta\ , \quad {\rm with} \,\,a\in\zb,\,\beta\in\mb\ ,
\end{equation} 
so the summation conditions listed in~\eqref{eq:5.8} read as follows:
\begin{equation}\label{eq:5.11} 
a+ \beta \in {\cal D},\,\,\, N(a +\beta)\le x\ ,
\end{equation} 
\begin{equation}\label{eq:5.12} 
a + \beta\equiv \rho\, (\mod F), \,\,\, a +\beta \equiv 0\,(\mod \fm )\ .
\end{equation} 
The congruence conditions show that not every $\beta \in \mb$ is 
admissible because the two residue classes $\beta- \rho\, (\mod F)$ 
and $\beta\,(\mod \fm )$ must both be represented by rational integers. 
Since $(F,\fm )=1$ we can choose one rational integer $\tilde a$ 
such that
\begin{equation}\label{eq:5.13} 
\tilde a +\beta\equiv \rho\, (\mod F), \quad \tilde a +\beta\equiv 0 
\,(\mod \fm )\ ,
\end{equation} 
and we get a single congruence condition for the variable $a$:  
\begin{equation}\label{eq:5.14} 
a\equiv \tilde a \,(\mod F\fm )\ .
\end{equation} 
 We shall not 
exploit these properties in any substantial fashion, but rather 
keep them in mind for consistency in the forthcoming arguments. 
Moreover, it is also advisable to keep in mind that all the conjugates 
of $\beta$ satisfy
\begin{equation}\label{eq:5.15} 
\beta^{(1)}, \ldots, \beta^{(n)} \ll x^{\frac{1}{n} }\ , 
\end{equation} 
which follows from~\eqref{eq:5.11}. 

From now on we think of $a$ as a variable which satisfies the 
conditions~\eqref{eq:5.11} and the congruence~\eqref{eq:5.14} 
while $\beta$ is inactive. 
Therefore, the conditions for $a$ which will emerge from the 
forthcoming transformations must be articulated rather precisely, but 
we do not need to be very explicit about the resulting features which 
depend only on $\beta$. For example, the signs $\pm$ which will come 
out of various applications of the reciprocity law will be independent of $a$ 
because $a$ runs over the fixed residue class~\eqref{eq:5.14}.

For $\fa= (\alpha)$ satisfiying~\eqref{eq:5.10}--\eqref{eq:5.14} we write 
$$
spin (\fa) = \Bigl(\frac{\alpha}{\alpha'}\Bigr) =
\Bigl(\frac{a + \beta}{a+\beta'}\Bigr) 
= \Bigl(\frac{\beta-\beta'}{a+\beta'}\Bigr)
$$
by the periodicity. If $\beta=\beta'$ we get no contribution, so we can 
assume
\begin{equation}\label{eq:5.16} 
\beta\neq \beta'\ . 
\end{equation} 
Next, we are going to interchange the upper entry $\beta - \beta'$ 
and the lower entry $a+\beta'$ by the reciprocity law. To do so, we first 
pull out from the ideal $(\beta -\beta')$ all its prime factors in $F$. 
We can write
\begin{equation}\label{eq:5.17} 
(\beta - \beta') =(\eta^2)\fc_0\fc_1\fc \ , 
\end{equation} 
with $(\fc,F)=1, \, \fc_1| F^{\infty}, \, (\fc_1,2)=1,\, \fc_0 |4^h$, 
where, as before, $h$ is the class number, 
$\eta\in {\cal O}, \, \eta | 2^{\infty}$. (This factorization is not unique 
but any such choice will lead to the same result for $A(x; \rho)$.) 

Let $\fC_0,\,\fC_1,\,\fC$ be three of the fixed ideals in ${\cal C}\ell$ which 
represent the inverse classes of $\fc_0,\,\fc_1,\,\fc$, respectively, 
see~\eqref{eq:2.8}--\eqref{eq:2.9}.
Therefore, $\fc_0 \fC_0 =(\gamma_0), \,\fc_1 \fC_1 =(\gamma_1)$ 
and $\fc \fC =(\gamma)$. 
Since $\fc_0 \fc_1\fc$ is a principal ideal, so is $\fC_0\fC_1\fC 
= (\gamma_2)$, say. Choosing appropriate associates we have
\begin{equation}\label{eq:5.18} 
\beta - \beta' =\eta^2\,\gamma_0\,\gamma_1\,\gamma\,\gamma_2^{-1}\ .
\end{equation} 
Hence,
$$
\Bigl(\frac{\beta-\beta'}{a+\beta'}\Bigr)
=\Bigl(\frac{\gamma_0\gamma_1\gamma_2\gamma}{a+\beta'}\Bigr)
=\Bigl(\frac{\gamma_0}{a+\beta'}\Bigr)\, 
\Bigl(\frac{\gamma_1\gamma_2\gamma}{a+\beta'}\Bigr)\ .
$$
Since $(\gamma_0)=\fc_0\fC_0$ divides $F/8$ and $a$ runs over a fixed residue 
class modulo $F$, the first symbol $(\gamma_0/(a+\beta'))$ does not depend 
on $a$; see Corollary~\ref{cor 2.1.2}. 
For the second symbol we apply the reciprocity law and obtain
$$
\Bigl(\frac{\gamma_1\gamma_2\gamma}{a+\beta'}\Bigr) =\pm 
\Bigl(\frac{a+\beta'}{\gamma_1\gamma_2\gamma}\Bigr)\ ,
$$
where, for fixed $\beta$, the sign $\pm$ does not depend on $a$. 
(This is the case because $a+\beta'$ is in a fixed odd class modulo 8, 
which suffices to imply that the relevant Hilbert symbols at even places 
do not depend on $a$. Moreover, at the infinite places the Hilbert symbols 
are equal to 1 because, by~\eqref{eq:5.11}, $a+\beta$ is totally positive.) 
Furthermore, we have
$$
\Bigl(\frac{a+\beta'}{\gamma_1\gamma_2\gamma}\Bigr) =
\Bigl(\frac{a+\beta'}{\fc}\Bigr)\, 
\Bigl(\frac{a+\beta'}{\fC_0\fc_1}\Bigr)
= \Bigl(\frac{a+\beta'}{\fc}\Bigr)\, 
\Bigl(\frac{\tilde a +\beta'}{\fC_0\fc_1}\Bigr)\ ,
$$
because $a\equiv \tilde a \, (\mod F)$ and $\fC_0\fc_1 | F^{\infty}$. 
We conclude that
\begin{equation}\label{eq:5.19} 
spin (\fa) = \pm \Bigl(\frac{a+\beta'}{\fc}\Bigr)\, 
\end{equation} 
where the $\pm$ sign does not depend on $a$ and where $\fc$ denotes 
the part of the ideal $(\beta -\beta')$ free of divisors of $F$, that is 
\begin{equation}\label{eq:5.20} 
\fc = \fc(\beta) = (\beta-\beta')/(\beta-\beta', F^{\infty})\ .
\end{equation} 

Hence, $A(x;\rho)$ splits as follows:
\begin{equation}\label{eq:5.21} 
A(x;\rho)= \sum_{\beta\in \mb}\pm T(x;\beta)\ , 
\end{equation} 
where $T(x;\beta)$ is given by
\begin{equation}\label{eq:5.22} 
T(x;\beta)= \sumflat_a\Bigl(\frac{a+\beta'}{\fc}\Bigr)\ . 
\end{equation} 
Here, the symbol $\flat$ means that $a$ runs over the rational integers 
satisfying the conditions~\eqref{eq:5.11} and the congruence~\eqref{eq:5.14}. 
Of course, these conditions impose some restrictions on $\beta\in\mb$. 
For example,~\eqref{eq:5.13} and~\eqref{eq:5.15} must hold or else the 
summation~\eqref{eq:5.22} is void. At this point we do not need to be 
very specific about the exact conditions for $\beta$.

We proceed to the estimation of $T(x;\beta)$. Our intention is to replace 
$\beta'$ in the upper entry of the symbol by a rational integer 
modulo $\fc$. This however may not be possible if the ideal $\fc$ contains 
prime divisors of degree greater than one. For this reason, we factor $\fc$ 
into
\begin{equation}\label{eq:5.23} 
\fc=\fg\fq\ , 
\end{equation} 
where $\fg$ takes from $\fc$ all prime ideals of degree greater than one, 
all ramified primes, and all unramified primes of degree one for which 
some different conjugate is also a factor of $\fc$. For $\fq$ taking the rest, 
note that $q = N\fq$ is a square-free number and norm $g=N\fg$ is a 
squarefull number co-prime with $q$. Let $b$ be a rational 
integer with $ b \equiv \beta'\, (\mod q)$. 
 This exists because  $q=N\fq$ is square-free. 
Note that $b$ is a rational integer 
which depends on $\beta$ but not on $a$. We have 
$$
\Bigl(\frac{a+\beta'}{\fc}\Bigr) =
\Bigl(\frac{a+\beta'}{\fg}\Bigr) \, \Bigl(\frac{a+\beta'}{\fq}\Bigr) 
= \Bigl(\frac{a+\beta'}{\fg}\Bigr) \, \Bigl(\frac{a+b}{\fq}\Bigr)\ . 
$$

Let $g_0$ be the product of all distinct prime divisors of $g$, 
\begin{equation}\label{eq:5.24} 
g_0=\prod_{p\mid g}p\ .
\end{equation} 
The quadratic residue symbol $(\alpha/\fg)$ is periodic in $\alpha$ 
modulo $\fg^* =\prod_{\fp|\fg}\fp$, hence it is periodic in $\alpha$ 
modulo $g_0$ because $\fg^*$ divides $g_0$. 
Therefore, the symbol $((a+\beta')/\fg)$ as a function 
of $a$ is periodic of period $g_0$. Splitting the sum~\eqref{eq:5.22} 
into residue classes modulo $g_0$, we get 
\begin{equation}\label{eq:5.25} 
T(x;\beta) \le \sum_{a_0\, (\mod g_0)} 
\Big|\sumflat_{a\equiv a_0\,(\mod g_0)}\Bigl(\frac{a+ b}{\fq}\Bigr)\Big|\ . 
\end{equation} 
Recall that the superscript $\flat$ indicates that the summation 
variable $a$ satisfies the conditions~\eqref{eq:5.11},~\eqref{eq:5.14}. 
These conditions imply
\begin{equation}\label{eq:5.26} 
a\ll x^{\frac 1n}\ ,
\end{equation} 
but we have to describe these conditions much more precisely. 
The first condition $\alpha = a +\beta\in {\cal D}$ is described 
by a system of linear inequalities $u\cdot \alpha > e\cdot \alpha$ 
for every $u\in 
\tilde{\cal U}$, where $\tilde{\cal U}$ is a finite subset of 
${\cal U}^+$. Hence, this condition means that $a$ runs over a single 
open interval whose endpoints depend on $\beta$. Next, we observe that the 
polynomial $(X+\beta^{(1)})\ldots (X+\beta^{(n)})$ has real coefficients,  
so the second condition $N(a+\beta) \le x$ means that $a$ runs over a 
collection of $n$ segments whose endpoints depend on $\beta$ and $x$. 
Therefore, the simultaneous conditions in~\eqref{eq:5.11} can be expressed 
by saying that $a$ runs over a certain collection of $n$ intervals.  

Finally, the congruence~\eqref{eq:5.14}, together with
$a\equiv a_0\,(\mod g_0)$, means that $a$ runs over a certain arithmetic 
progression of modulus $k$ which divides $g_0F N\fm$. We can 
assume that $m=N\fm$ and 
$q=N\fq$ are co-prime. If not, there is a prime ideal $\fp$ such that 
$\fp|\fm'$ and $\fp|q$. Let $\fp^{\tau}$ be a conjugate of $\fp$ which 
divides the ideal $\fq$. Then, the quadratic residue symbol 
in~\eqref{eq:5.25} factors as follows:
$$
\Bigl(\frac{a+ b}{\fq}\Bigr)= \Bigl(\frac{a+ b}{\fp^{\tau}}\Bigr)\ldots 
= \Bigl(\frac{a+ b}{\fp}\Bigr)\ldots = \Bigl(\frac{a+ \beta'}{\fp}\Bigr)\ldots 
= 0 
$$
because $b\equiv \beta'\, (\mod q)$ and $\tilde a + \beta' \equiv 0\, 
(\mod \fm')$. Since $g_0, \, F$ are also co-prime with $q$, we have 
$(g_0Fm,q)=1$ so $a$ in~\eqref{eq:5.25} runs over an arithmetic progression 
of modulus $k$ which is co-prime with $q$.  

Having said these things, we see that the inner sum 
in~\eqref{eq:5.25} can be arranged as the sum of $n$ sums, each of which 
runs over a single segment of an arithmetic progression 
of length $\ll x^{1/n}$. Since 
$\chi_{\fq}(\ell)=(\ell/\fq)$ is a real Dirichlet character of modulus 
$q=N\fq$ (see~\eqref{eq:2.15} and Lemma~\ref{lem 2.6}) we have here $n$ 
incomplete character sums of length $\ll x^{1/n}$ and the modulus $q\ll x$ 
of the character is co-prime to the modulus $k$ of the progression. 
Therefore, if $q$ is not squarefull, Conjecture $C_n$ applies, 
(or rather its consequence for arithmetic progressions as described 
at the end of Section~\ref{sec 8}), giving (see~\eqref{eq:8.4}),  
\begin{equation}\label{eq:5.26.1} 
T(x;\beta) \ll g_0x^{\frac{1-\delta}{n}+\ve}\ . 
\end{equation} 
Here, the implied constant depends only on $\ve$ and the field $K$ 
but not on $\beta$. Perhaps the reader is wondering why the implied constant 
in~\eqref{eq:5.26.1} does not depend on $\beta$, although many steps before 
we arrived there depend on $\beta$. The answer is that the 
estimate~\eqref{eq:8.4} in Conjecture $C_n$ holds for any incomplete 
character sum, regardless of the location of the segment where the summation 
takes place. This feature is really vital for our application because
in~\eqref{eq:5.25} we have no idea where $b$, the rational representative 
of $\beta'$ modulo $\fq$, can be in relation to $q=N\fq$. 

Recall that~\eqref{eq:5.26.1} holds provided $q=N\fq$ is not squarefull. 
However, if $q$ is squarefull, then $q=1$ and 
\begin{equation}\label{eq:5.26.2} 
F^2N(\beta -\beta')\,\,\, {\rm is \,\,squarefull}\ .
\end{equation} 
The condition~\eqref{eq:5.26.2} is satisfied very rarely, so we can count 
crudely. Denote by $A_{\square}(x;\rho)$ the contribution to 
$A(x;\rho)$ of the terms for which~\eqref{eq:5.26.2} holds. We can afford 
to ignore the congruence conditions~\eqref{eq:5.12}. We have 
$$
A_{\square}(x;\rho)\le |\{\alpha\in{\cal D}; \,  N\alpha \le x,\, 
F^2N(\beta -\beta')\,\, {\rm squarefull}\,\}|\ . 
$$

Since $\alpha \in {\cal D}, \, N(\alpha) \le x$ all the conjugates 
$\alpha^{(k)}$ are $\ll x^{1/n}$. Hence, $ |a|\le y$ and all the conjugates of 
$\gamma =\beta - \beta'=\alpha -\alpha'$ satisfy $|\gamma^{(k)}| \le y$ 
for some $y\asymp x^{1/n}$. Recall that the map $\mb\rightarrow {\cal O}$ 
given by $\beta \rightarrow \beta -\beta'$ is injective (see 
Lemma~\ref{lem 2.4}). Therefore, we have 
$$
A_{\square}(x;\rho)\le y|\{\gamma\in{\cal O}; \,  |\gamma^{(k)}| \le y, \, 
F^2N(\gamma)\,\,\, {\rm squarefull}\,\}|\ . 
$$ 
Here, we can replace the counting of integers $\gamma$ by the counting 
of the principal ideals they generate, each integer occurring with 
multiplicity $\ll (\log x)^n$ by Lemma~\ref{lem 2.3}. Hence 
$$
A_{\square}(x;\rho)\ll x^{\frac 1n}(\log x)^n 
|\{\fb\subset{\cal O}; \,  N\fb \le X, \, 
F^2N\fb \,\,\, {\rm squarefull}\,\}|\ , 
$$ 
where $X=y^n$ so $X\asymp x$. Note that we have moved from the integers 
in the submodule $\lb\subset {\cal O}$ (see Section~\ref{sec 2})  
to the ideals in ${\cal O}$ 
because detecting squarefull norms in a submodule of lower rank could 
be very difficult. Fortunately, we can afford the loss which results from 
this extension because $n\ge 3$. Now, we can exploit the multiplicative 
structure of the ideals in ${\cal O}$ which gives us the bound
$$
A_{\square}(x;\rho)\ll x^{\frac 1n}(\log x)^n 
\sum_{\substack{ b \le X \\ 
F^2b\,\,\, {\rm squarefull}}}\tau_n(b)\ , 
$$
where $b$ runs over positive rational integers and $\tau_n(b)$ denotes 
the divisor function of degree $n$ so that $\tau_n(b)\ll b^{\ve}$. 
Hence, we conclude 
\begin{equation}\label{eq:5.27} 
A_{\square}(x;\rho)\ll x^{\frac 12 +\frac 1n +\ve}\ ,
\end{equation} 
where the implied constant depends on $\ve$ and the field $K$. 

Let $A_0(x;\rho)$ be the contribution to $A(x;\rho)$ of the terms
$\alpha = a+\beta$ for which~\eqref{eq:5.26.2} does not hold. Therefore, 
we have the following partition: 
\begin{equation}\label{eq:5.28} 
A(x;\rho)= A_{\square}(x;\rho) + A_0(x;\rho)
\end{equation} 

To estimate $A_0(x;\rho)$ we can use~\eqref{eq:5.26.1} for every relevant 
$\beta$, but the bound~\eqref{eq:5.26.1} is useless for $g_0$ too large. 
Therefore, we make the further partition 
\begin{equation}\label{eq:5.29} 
A_0(x;\rho)= A_1(x;\rho) + A_2(x;\rho)+ A_3(x;\rho)\ ,
\end{equation} 
where the components run over $\alpha = a+\beta,\, \beta\in\mb$ with 
$\beta$ such that 
\begin{equation}\label{eq:5.30} 
g_0\le Z \quad {\rm in}\,\, A_1(x;\rho)\ ,
\end{equation} 
\begin{equation}\label{eq:5.31} 
g_0 > Z, \,\, g\le Y  \quad {\rm in}\,\, A_2(x;\rho)\ ,
\end{equation} 
\begin{equation}\label{eq:5.32}
g_0 > Z, \,\, g> Y  \quad {\rm in}\,\, A_3(x;\rho)\ .
\end{equation} 
We shall choose $Z\le Y$ later. 

To estimate $A_1(x;\rho)$ we use~\eqref{eq:5.26.1} and then sum 
over $\beta\in\mb$ satisfying~\eqref{eq:5.15}, ignoring the other 
restrictions. This gives 
\begin{equation}\label{eq:5.33} 
A_1(x;\rho)\ll Zx^{1-\frac{\delta}{n}+\ve}\ .
\end{equation} 

The estimation of $A_3(x;\rho)$ is also very quick. We treat $A_3(x;\rho)$
by arguments similar to those we applied to $A_{\square}(x;\rho)$. The 
condition that $F^2N(\gamma)$ is squarefull is 
now replaced by $g|N(\gamma)$.  
Hence, 
\begin{equation*}
\begin{split}
A_3(x;\rho) & \ll x^{\frac 1n}\,(\log x)^n\sum_{\substack{g>Y\\ 
g\,\,{\rm squarefull}}}\sum_{\substack{b\le X\\g|b}}\tau_n(b) \\ 
& \ll x^{1+\frac 1n +\ve}
\sum_{\substack{Y<g\le x\\g\,\,{\rm squarefull}}}g^{-1}\ . 
\end{split}
\end{equation*}
This last sum is estimated by 
$$
Y^{-\frac 12 }\sum_{\substack{g\le x\\g\,\,{\rm squarefull}}}g^{-\frac 12}
\le Y^{-\frac 12 }\prod_{p\le x}
\Bigl(1+\frac 1p \, \bigl(1-\frac{1}{\sqrt p}\bigr)^{-1}\Bigr) 
\ll Y^{-\frac 12 } \log x\ .
$$
Hence, we conclude that
\begin{equation}\label{eq:5.34}
A_3(x;\rho)\ll Y^{-\frac 12}x^{1+\frac 1n +\ve}\ .
\end{equation} 

It remains to estimate $A_2(x;\rho)$ which is quite a difficult job. 
If, along the lines of the proof of~\eqref{eq:5.34}, we had not 
wasted the information that $\beta-\beta'$ is in a submodule 
$\lb\subset {\cal O}$ of rank $n-1$, then the factor $x^{1/n}$ 
in~\eqref{eq:5.34} would be saved, so the result would be useful for 
$Y>x^{\delta}$ and there would be no need to consider the sum $A_2(x;\rho)$
over the middle range~\eqref{eq:5.31}. But, it is difficult to take 
this information into account when $g=N\fg$ is quite large. We are now 
going to exploit this information for the estimation of $A_2(x;\rho)$ 
but the arguments need to be more delicate. 

We begin with the arguments that were applied to $A_3(x;\rho)$ but 
keep the condition $\gamma =\beta-\beta' \equiv 0\,(\mod \fg)$ in place of 
$N(\gamma)\equiv 0\, (\mod g)$. We obtain
\begin{equation}\label{eq:5.35}
|A_2(x;\rho)|\le y\sum_{\substack{\fg\\g_0>Z,\,g\le Y}}E_{\fg}(y)\ ,
\end{equation} 
where $y\asymp x^{1/n}$ and 
$$
E_{\fg}(y)=|\{\gamma\in\lb;\,\gamma\equiv 0\,(\mod \fg),\,
|\gamma^{(k)}|\le y\,\, {\rm for\,\,all} \,\, k\}|\ .
$$
Recall that $\fg$ (see~\eqref{eq:5.23}) runs over ideals, 
all of whose prime factors have degree greater than one or, if of degree 
one, have another prime factor of the same norm, that 
$g=N\fg$ and $g_0$ is the product of all the 
distinct primes in $g$. Note that every prime ideal dividing $\fg$ 
also divides $g_0$. 

For the estimation of $E_{\fg}(y)$ we express every $\gamma$ in terms 
of the basis $\eta_2,\ldots,\eta_n$ of $\lb$, getting 
$\gamma= a_2\eta_2 +\ldots +a_n\eta_n$ with integer coefficients 
$a_2,\ldots,a_n$ satisfying
\begin{equation}\label{eq:5.36}
a_2,\ldots,a_n \ll y\ ,
\end{equation} 
\begin{equation}\label{eq:5.37}
a_2\eta_2 +\ldots +a_n\eta_n \equiv 0\,(\mod \fg)\ . 
\end{equation} 

Next, we split the coefficients $a_2,\ldots,a_n$ according to their 
residue classes modulo $g_0$, say $r_2,\ldots,r_n$. For each 
unramified prime $\fp |\fg$ of degree greater than one we take an automorphism 
$\tau \neq\,\, {\rm id}$ such that $\fp^{\tau} 
=\fp$. Such a $\tau$ exists (in the decomposition 
group of the prime $\fp$), because deg $\fp =f > 1$. 
Then,~\eqref{eq:5.37} yields two different congruences modulo $\fp$:
\begin{equation}\label{eq:5.37.1}
\begin{split}
r_2\eta_2 +\ldots +r_n\eta_n & \equiv 0\,(\mod \fp)\ , \\
r_2\eta_2^{\tau} +\ldots +r_n\eta_n^{\tau} & \equiv 0\,(\mod \fp)\ .
\end{split}
\end{equation}
For each unramified prime $\fp$ dividing $\fg$ of degree one there exists 
a companion $\fp^{\tau^{-1}} \neq \fp$ which also divides $\fg$. 
Then,~\eqref{eq:5.37} yields one congruence modulo $\fp\, \fp^{\tau^{-1}}$: 
$$
r_2\eta_2 +\ldots +r_n\eta_n  \equiv 0\,(\mod \fp\,\fp^{\tau^{-1}})\ .
$$
This single congruence yields the same system~\eqref{eq:5.37.1}. 

By Lemma~\ref{lem 2.5.1} we have $\eta_k\eta_{\ell}^{\tau}\neq 
\eta_k^{\tau}\eta_{\ell}$ for some $2\le k\neq\ell\le n$. Hence, if 
$\fp$ does not divide $\fh$ (see the definition of the ideal $\fh$ 
in~\eqref{eq:2.14.7})
$$
\eta_k\eta_{\ell}^{\tau}-\eta_k^{\tau}\eta_{\ell}\not\equiv 0\,(\mod \fp)\ ,
$$
so $r_k,\,r_{\ell}$ are determined uniquely modulo $\fp$ by the other 
residue classes. But $r_k,\,r_{\ell}$ are rational so they are determined 
uniquely modulo every prime ideal which is a conjugate of $\fp$. Hence, 
by the Chinese Remainder Theorem, $r_k,\,r_{\ell}$ are determined 
uniquely modulo $p$ for every $p|g_0$, $p$ co-prime with $h=N\fh$, 
because $p$ is the product of distinct conjugates of $\fp$. 
Thus,
$$
|\{r_2,\ldots,r_n \, (\mod p)\,\, {\rm admissible}\}| = p^{n-3}\ .
$$
For the primes $p|h$ we use the trivial bound $p^{n-1}$, so the total 
number of admissible residue classes satisfies 
\begin{equation}\label{eq:5.38}
|\{r_2,\ldots,r_n \, (\mod g_0)\,\, {\rm admissible}\}| \ll g_0^{n-3}\ ,
\end{equation}
where the implied constant depends only on the field $K$. 
Hence, the number of coefficients $a_2,\ldots, a_n$ 
satisfying~\eqref{eq:5.36} and~\eqref{eq:5.37} is bounded by 
$$
g_0^{n-3}\bigl(yg_0^{-1} + 1\bigr)^{n-1}\ll y^{n-1} g_0^{-2} + g_0^{n-3}\ ,
$$
and so
\begin{equation}\label{eq:5.39}
E_{\fg}(y)\ll y^{n-1}g_0^{-2} + g_0^{n-3}\ .
\end{equation} 
Introducing~\eqref{eq:5.39} into~\eqref{eq:5.35} we obtain
$$
|A_2(x;\rho)|\ll\sum_{\substack{\fg\\g_0>Z,\,g\le Y}}
\bigl(y^{n}g_0^{-2} + yg_0^{n-3}\bigr)\ .
$$
The number of ideals $\fg$ with $N\fg = g$ is bounded by 
$\tau(g)\ll Y^{\ve}$. Hence,
$$
|A_2(x;\rho)|\ll x^{\ve}
\sum_{\substack{g\,\, {\rm squarefull}\\ g_0>Z,\,\,g\le Y}}
\bigl(xg_0^{-2} + x^{\frac 1n}g_0^{n-3}\bigr)\ , 
$$
where $g_0$ is the product of all distinct prime divisors of $g$. 

We estimate the first sum as follows:
$$
\sum_g g_0^{-2} \le Z^{-1}\sum_g g_0^{-1}
\le Z^{-1}\prod\Bigl(1+\frac 1p \, \bigl(1-\frac 1p \bigr)^{-1}\Bigr) 
\ll Z^{-1}\log Y\ . 
$$
In the second sum we use the fact that $g_0^2|g$, so 
$g_0\le \sqrt g \le \sqrt Y$, whence
$$
\sum_g g_0^{n-3} \le Y^{\frac{n-3}{2}}
\sum_{\substack{g\,\, {\rm squarefull}\\ g\le Y}} 1 \ll Y^{\frac n2 -1}\ . 
$$
We conclude that
\begin{equation}\label{eq:5.40}
A_2(x;\rho) \ll Z^{-1}x^{1+ \ve}+ Y^{\frac n2 -1}x^{\frac 1n + \ve}\ ,
\end{equation} 
where the implied constant depends on $\ve$ and the field $K$. 

Inserting the three 
estimates,~\eqref{eq:5.33},~\eqref{eq:5.34},~\eqref{eq:5.40}, 
into~\eqref{eq:5.29}, we deduce that  
$$
A_0(x;\rho) \ll x^{\ve}\,\bigl(Zx^{1-\frac{\delta}{n}} +Y^{-\frac 12} 
x^{1+\frac 1n} +Z^{-1}x +Y^{\frac n2 -1}x^{\frac 1n}\bigr)\ , 
$$
where $Z\le Y$ are at our disposal. We choose $Z= x^{\delta/2n}$ and 
$Y=x^{2/(n-1)}$, getting 
\begin{equation}\label{eq:5.41}
A_0(x;\rho) \ll x^{1- \frac{\delta}{2n}+ \ve}\ . 
\end{equation} 
Adding to~\eqref{eq:5.41} the bound~\eqref{eq:5.27}, we find 
that~\eqref{eq:5.41} holds also for $A(x;\rho)$. Finally, summing this 
bound over the residue 
classes $\rho \,(\mod F)$ and adding the bound~\eqref{eq:5.9}, 
we complete the proof of the bound~\eqref{eq:5.3}, that is 
Proposition~\ref{prop 5.1}.

\section{Bilinear Forms with Spin}\label{sec 6}

Let $F$ be a fixed positive integer which is a multiple of $\fM$ and $\frak f$ 
(see~\eqref{eq:2.9}. We consider a general bilinear form
\begin{equation}\label{eq:6.1}
{\cal B}(x,y)=\ssum_{\substack{
(\frak a\frak b,F)=1 \\ N\frak a \le x,\, N\frak b\le y}}
v_{\frak a}\,w_{\frak b}\, r(\fa\fb)\, spin(\frak a\, \frak b)\ ,
\end{equation}
where $v_{\frak a}, \, w_{\frak b}$ are arbitrary complex numbers with 
$|v_{\frak a}| \le 1, \, |w_{\frak b}| \le 1 $. Recall that the factor
$r(\frak a\frak b)$ means that $\frak a\frak b$ is a principal ideal 
which has a totally positive generator in the residue class $\mu$ 
modulo $\fM$. 

\begin{prop}\label{prop 6.1}
Let $K/ \qb$ be a totally real Galois extension of degree $n\ge 3$. Then,
\begin{equation}\label{eq:6.2}
{\cal B}(x,y) \ll (x+y)^{\frac{1}{6n}}(xy)^{1- \frac{1}{6n} +\ve}
\end{equation}
for any $x,\, y\ge 2$ and $\ve >0$, the implied constant depending 
only on $\ve$ and the field $K$.
\end{prop}

\begin{proof}
Since the coefficients $v_{\frak a}, \, w_{\frak b}$ are arbitrary, we 
can assume without loss of generality that they are supported on fixed 
ideal classes, say 
\begin{equation}\label{eq:6.4}
\frak a\frak A =(\alpha)\ , \quad \alpha \succ 0\ ,
\end{equation}
\begin{equation}\label{eq:6.5}
\frak b \frak B = (\beta)\ , \quad \beta \succ 0\ ,
\end{equation}
as in~\eqref{eq:2.9.1}, with $\fA\neq\fB$. 
Furthermore we can assume that $\alpha, \, \beta$ have fixed residue classes 
modulo 8. Then, $spin\, (\frak a \frak b)$ factors as in~\eqref{eq:2.10}. 
The middle symbols in~\eqref{eq:2.10} separate $\frak a$ from $\frak b$ 
so they can be attached to the coefficients $v_{\frak a}, \, w_{\frak b}$. 
Hence we have
$$
|{\cal B}(x,y)| \le \sum_{\substack{(\frak a , F)=1\\ N\frak a \le x} }
\Big| \sum_{\substack{(\frak b, F)=1\\ N\frak b \le y}} w_{\frak b}
\Bigl(\frac {\alpha}{\frak b' \frak b^-}\Bigr)\Big|\ .
$$  
Here, the coefficients $v_{\fa}$ have disappeared and
 $w_{\frak b}$ are complex numbers with $|w_{\frak b}| \le 1$, not 
necessarily the original ones in~\eqref{eq:6.1}. Note that we camouflaged 
the fact that $\frak b$ satisfies~\eqref{eq:6.5} by incorporating this 
information into the coefficients $w_{\frak b}$ because we no longer will 
make any use of this. Note also that the characteristic function 
$r(\fa\, \fb)$ has disappeared for the same reason.  
Moreover, we do not need all the information about the ideal 
$\frak a$ coming from~\eqref{eq:6.4} but we must keep~\eqref{eq:6.4} in 
mind for a while until the relation of 
$\alpha$ being a function of $\frak a$ disappears.
  
Our strategy is to create another bilinear form, one in which the variables 
have very different sizes. To this end we apply a high power 
H\"older inequality and use the multiplicativity of the quadratic residue 
symbol with respect to the modulus, which holds by definition. We obtain 
$$
|{\cal B}(x,y)|^k \le {\cal N}(x)^{k-1}
\sum_{\substack{(\frak a,F)=1\\ N\frak a \le x}}
\Big|\sum_{\substack{(\frak b,F)=1\\ N\frak b \le y}}w_{\frak b} 
\Bigl(\frac{\alpha}{\frak b' \frak b^-}\Bigr)\Big|^k\ , 
$$
where ${\cal N}(x)$ is the number of ideals ${\frak a} $ 
with $N{\frak a}\le x$, so ${\cal N}(x)\ll x$. Here, $k$ is a 
positive integer to be chosen later, not necessarily even. Next, we get 
$$
\sum_{\frak a}\big| \sum_{\frak b}\big|^k \le 
\sum_{(\frak c,F)=1}\Big| \sum_{\substack{(\frak a,F)=1\\ N\frak a \le x}}
\ve (\frak a)\Bigl(\frac{\alpha}{\frak c' \frak c^-}\Bigr)\Big|\ , 
$$   
where $|\ve (\fa)|=1$
(precisely, $\bar{\ve} (\frak a)$ is the $k$-th power of the complex sign 
of the inner sum over $\fb$),
and $\frak c= \frak b_1\ldots\frak b_k$ 
with $\frak b_1,\, \ldots, \, \frak b_k$ running independently over 
all ideals of norm $\le y$. Therefore, $N\frak c\le y^k=Y$, say, and 
the number of representations of $\frak c$ as the product of $k$ ideals is 
$\tau_k(\frak c) \ll y^{\ve}$. Hence we find 
$$
|{\cal B}(x,y)|^k \ll y^{\ve}x^{k-1}
\sum_{\substack{(\frak c,F)=1\\N\frak c \le Y}}
\Big| \sum_{\substack{(\frak a,F)=1\\ N\frak a \le x}}
\ve (\frak a)\Bigl(\frac{\alpha}{\frak c' \frak c^-}\Bigr)\Big|\ , 
$$ 
where the implied constant depends on $\ve, \, k$ and the field $K$. 
At this point our mission is accomplished because, if $k$ is sufficiently 
large, then $Y=y^k$ is much larger than $x$. To take advantage of this 
disproportion of the sizes of $\frak c$ and $\frak a$, we intend to execute 
the summation over $\frak c$ first while holding $\frak a$ inactive. This 
requires however an interchange of the positions of $\frak c$ and $\frak a$ 
which can be accomplished by the reciprocity law, followed by an 
application of Cauchy's inequality.

For the use of the reciprocity law we must split $\frak c$ into ideal classes 
and, for each class, treat the corresponding sum separately. Let this class 
be determined by
\begin{equation}\label{eq:6.6}
\frak c\, \frak C =(\gamma)\ , \quad \gamma \succ 0\ ,
\end{equation}
where $\frak C \in {\cal C}\ell$ is the chosen ideal 
which is different from $\fA$ in~\eqref{eq:6.4}. Such a choice is possible 
because every ideal class has two representatives in ${\cal C}\ell$. Note 
that $N\fC$ is co-prime with $N\fA$ because of~\eqref{eq:2.9}. 
Actually, if we kept the information~\eqref{eq:6.5} throughout, then 
we would already know the ideal class of $\frak c$, namely the inverse 
class of $\frak B^k$, but this information would not save us much work. 
Thus, we appeal to~\eqref{eq:6.6} to obtain 
$$
\Bigl(\frac{\alpha}{\frak c' \frak c^-}\Bigr) =
\Bigl(\frac{\alpha}{\gamma' \gamma^-}\Bigr)\,
\Bigl(\frac{\alpha}{\frak C' \frak C^-}\Bigr) 
= \pm\, \Bigl(\frac{\gamma' \gamma^-}{\frak a}\Bigr)\, 
\Bigl(\frac{\alpha}{\frak C' \frak C^-}\Bigr)\ ,
$$
where the sign $\pm$ depends only on the residue classes of 
$\alpha, \, \gamma$ modulo 8. Hence,
$$
|{\cal B}(x,y)|^k \ll y^{\ve}x^{k-1}
\sum_{\substack{(\frak c,F)=1\\N\frak c \le Y}}
\Big| \sum_{\substack{(\frak a,F)=1\\ N\frak a \le x}}
\ve (\frak a)\Bigl(\frac{\gamma'\,\gamma^-}{\frak a}\Bigr)\Big|\ , 
$$ 
where $|\ve (\frak a)|=1$. Here, we did not display the sign $\pm$ 
because its dependence on $\gamma \, (\mod 8)$ disappears by 
the positivity argument in our summation and its dependence on 
$\alpha \, (\mod 8)$ is absorbed by the floating coefficient $\ve (\frak a)$. 

Now, we can forget the condition~\eqref{eq:6.4} for $\frak a$, which 
we had kept in mind until now, by building it into the coefficient 
$\ve (\frak a)$. We could not have done this earlier because $\alpha$, 
the generator of $\frak a \frak A$ was present in the quadratic residue 
symbol. 

If we require in~\eqref{eq:6.6}, as we may, that $\gamma$ belong to the 
fundamental domain ${\cal D}^*$, then there is a one-to-one correspondence 
between $\frak c$ and $\gamma$ with $\gamma \equiv 0\, (\mod \frak C)$. 
Hence,
\begin{equation}\label{eq:6.7}
|{\cal B}(x,y)|^k \ll y^{\ve}x^{k-1}{\cal E}(X, x)
\end{equation}
where 
\begin{equation}\label{eq:6.8}
{\cal E}(X, x) =
\sum_{\substack{\gamma \in \bar{\cal D}\\N\gamma\, \le X}}
\Big| \sum_{\substack{\frak a\,\,\, {\rm odd}\\ N\frak a\, \le x}}
\ve (\frak a)\Bigl(\frac{\gamma'\gamma^-}{\frak a}\Bigr)\Big|\ , 
\end{equation}
where $X=YF$. Note that we ignored the conditions
$(\gamma,F)=1, \, \gamma\equiv 0\, (\mod \frak C)$ and extended ${\cal D}^*$ 
to its closure $\bar{\cal D}$ as we may, by positivity. Have in mind that $F$ 
is a constant, depending on the field $K$, so that $X\asymp Y=y^k$. 

Applying Cauchy's inequality and changing the order of summation, we 
arrive at 
\begin{equation}\label{eq:6.9}
\begin{split}
|{\cal E}(X, x)|^2 
& \le \Bigl(\sum_{\substack{\gamma \in \bar{\cal D}\\N\gamma \le X}}\, 1 \Bigr)
\sum_{\substack{\fq\,\, {\rm odd}\\ N\frak q \le x^2}}\Big|
 \sum_{\fa_1\fa_2 =\fq}\ve (\frak a_1)\bar {\ve} (\frak a_2)
 \sum_{\substack{\gamma \in \bar{\cal D}\\N\gamma \le X}}
\Bigl(\frac{\gamma'\gamma^-}{\fq}\Bigr)\Big| \\
& \ll X\sum_{\substack{\fq\,\, {\rm odd}\\ N\frak q \le x^2}}
\tau(\frak q)|
\Sigma_{\frak q}(X)|\ .
\end{split} 
\end{equation}
Here, $\tau (\frak q)$ gives a bound for the number
of representations $\frak q =\frak a_1\frak a_2$ and 
\begin{equation}\label{eq:6.10}
\Sigma_{\frak q}(X) = \sum_{\substack{\gamma 
\in \bar{\cal D}\\N\gamma \,\le X}}
\Bigl(\frac{\gamma'\gamma^-}{\frak q}\Bigr)\ .  
\end{equation}
Let $q=N\frak q$ be the norm of $\frak q$. Note that the diagonal 
terms $\fa_1=\fa_2$ are covered by this case. If $q$ is squarefull, we apply 
the trivial bound
$$
\Sigma_{\frak q}(X) \ll X\ .
$$
If $q$ is not squarefull, that is $q$ has a prime factor $p$ 
whose square does not divide $q$, 
we split $\Sigma_{\frak q}(X)$ into residue classes modulo $\frak q$: 
$$
\Sigma_{\frak q}(X)= \sum_{\nu\, (\mod\, \frak q)}
\Bigl(\frac{\nu'\nu^-}{\frak q}\Bigr)\, {\cal N}
(x;\frak q, \nu)\ .
$$
By Lemma~\ref{lem 3.6} we obtain 
$$
\Sigma_{\frak q}(X)=\frac{{\cal N}(x)}{q} \sum_{\nu\, (\mod\, \frak q)}
\Bigl(\frac{\nu'\nu^-}{\frak q}\Bigr)
+O\bigl(q X^{1-\frac 1n }\bigr)\ .
$$
Here, the complete sum over $\nu\, (\mod \fq)$ is equal to 
$q^{-1}$ times the complete sum over $\alpha \, ( \mod \fq'\,\fq^{-})$ 
in~\eqref{eq:2.13}, so it vanishes and we are left with 
$$
\Sigma_{\frak q}(X)\ll q X^{1-\frac 1n }\ . 
$$
Hence, by~\eqref{eq:6.9} we derive 
\begin{equation*}
\begin{split}
|{\cal E}(X, x)|^2 & \ll X^2\,\sum_{\substack{N\frak q \le x^2\\ 
N\frak q\,\, {\rm squarefull}}} 1 + X^{2-\frac 1n }\sum_{N\frak q \le x^2} 
\tau(\frak q) q \\ 
& \ll y^{2k}\bigl(x+y^{-k/n}x^4\bigr)x^{\ve}\ , 
\end{split}
\end{equation*}
so that 
\begin{equation}\label{eq:6.11}
{\cal E}(X, x) \ll y^{k}\bigl(x^{1/2}+y^{-k/2n}x^2\bigr)x^{\ve}\ . 
\end{equation}
Inserting the result into~\eqref{eq:6.7}, we find 
that 
$$
{\cal B}(x,y) \ll (xy)^{1+\ve}\bigl(x^{-1/2k}+y^{-1/2n}x^{1/k}\bigr)\ .
$$ 
We choose $k=3n$, obtaining
\begin{equation*}
\begin{split}
{\cal B}(x,y) & \ll (xy)^{1+\ve}\bigl(x^{-1/6n}+y^{-1/2n}x^{1/3n}\bigr)\\ 
& < (xy)^{1- 1/6n +\ve}(x+y)^{1/6n}\bigl(1+ y^{-1}x\bigr)^{1/3n}\ .
\end{split}
\end{equation*}
By the symmetry of the bilinear form we can assume that $x\le y$, 
getting~\eqref{eq:6.2}. 
\end{proof}

\section{Conclusion of Proof}\label{sec 7}

We now have all the pieces for the application of Proposition~\ref{cor 4.2} 
to the sum of the spins of prime ideals. Specifically, we apply 
Proposition~\ref{cor 4.2} for the sequence ${\cal A} = (a_{\fn})$ 
with $a_{\fn}= 1$ if $(\fn, F)= 1, \, \fn = (\alpha), \, \alpha\succ 0, \, 
\alpha \equiv \mu \, (\mod\, \fM)$ and $a_{\fn} = 0$ otherwise. 

First,~\eqref{eq:5.3} gives 
us~\eqref{eq:4.7} with $\vartheta =\delta/2n$. Next,~\eqref{eq:6.2} 
gives us~\eqref{eq:4.8} with $\theta = 1/6n$. Therefore,~\eqref{eq:4.9} 
becomes 
\begin{equation}\label{eq:7.1} 
\sum_{\substack{N\frak n \le x\\ (\frak n,F)=1}}
\Lambda (\frak n )\,r(\fn)\,  spin(\frak n) \ll x^{1-\delta/2n(12n+1) + \ve}\ .
\end{equation}
Here, the condition $(\frak n,F)=1$ can be removed because the powers 
of prime ideals which divide $F$ contribute a negligible quantity. 
Moreover, the powers of primes which are not primes also 
contribute a negligible amount, so we are left with $\frak n =\frak p$, 
each one weighted by $\log N\frak p$. These logarithmic weights can be 
removed by partial summation. Hence~\eqref{eq:7.1} implies~\eqref{eq:1.2}.   

\section{Direct Estimates for Character Sums}\label{sec 8}

Our arguments in this paper are powered by estimates for real character 
sums over short intervals,
\begin{equation}\label{eq:8.1}
S_{\chi}(M,N) = \sum_{M< n \le M+N} \chi(n),\, \quad N\ge 1\ . 
\end{equation}
If $\chi\, (\mod q)$ is not principal, one expects to beat the trivial bound 
$S_{\chi}(M,N) \ll N$ by exploiting some cancellation due to the random 
sign changes of $\chi (n)$. The celebrated result of D. Burgess [Bu] gives us 
the bound
\begin{equation}\label{eq:8.2}
S_{\chi}(M,N)  \ll N^{1-\frac 1r}q^{\frac{r+1}{4r^2} +\ve} ,
\end{equation}
with any integer $r\ge 1$ and any $\ve >0$, the implied constant depending 
only on $r$ and $\ve$. The bound becomes trivial if $N\le q^{(r+1)/4r}$ so, 
no matter how large we choose $r$ in~\eqref{eq:8.2}, we get nothing useful 
for sums of length $N\le q^{1/4}$. Fortunately, for our application to 
estimate the sums of spins in this paper (see Section~\ref{sec 5}), 
in the case of a 
cubic field we encounter sums of length $N$ as large as $q^{1/3}$ and we can 
appeal to Burgess' estimate. 
We take~\eqref{eq:8.2} with $r=6$ to obtain:
\begin{cor}\label{cor 8.1}
Let $\chi\,(\mod q)$ be a non-principal real character. Then
\begin{equation}\label{eq:8.3}
S_{\chi}(M,N)  \ll N^{\frac 56}q^{\frac{7}{144} +\ve} ,
\end{equation}
with any $\ve>0$, the implied constant depending only on $\ve$. 
\end{cor}

When the degree of the field $K$ is $n>3$ we need a non-trivial bound  
for $S_{\chi}(M,N)$ with $N$ of size about $q^{1/n}$ and nothing useful 
is available at present for such short sums (except in the case of 
special moduli, cf [IK]). Thus, in order to cover the fields of higher degree 
we have no option other than to postulate adequate estimates for these 
short character sums. 

\noindent
{\bf Conjecture}\, $C_n$: Let $n\ge 3, \, Q \ge 3,\, N\le Q^{1/n}$. 
For any real 
non-principal character $\chi\, (\mod q)$ of modulus $q\le Q$ we have 
\begin{equation}\label{eq:8.4}
S_{\chi}(M,N)  \ll Q^{\frac{1-\delta}{n} +\ve} ,
\end{equation} 
with some $\delta=\delta (n)>0$ and any $\ve >0$ the implied constant 
depending only on $\ve$ and $n$.

\noindent
{\bf Remarks}: 
\, Burgess stated his estimate~\eqref{eq:8.2} for any $r\ge 1$ and 
any non-principal character $\chi\, (\mod q)$, but, in the case of 
general $r$, only for $q$ cube-free. However, if $\chi$ is real, the character 
sum does not change if $q$ is reduced, as it always can be, 
to some cube-free divisor of $q$ by dividing out the largest square which 
does not remove any prime factor completely. Thus~\eqref{eq:8.3} is correct 
as stated. Actually, we are only using these bounds in the case of 
square-free modulus.

If $N\le Q^{1/3}$ and $q\le Q$ then~\eqref{eq:8.3} implies~\eqref{eq:8.4} 
proving Conjecture $C_3$ with 
\begin{equation}\label{eq:8.6}
\delta = \delta (3) =\frac{1}{48}\ . 
\end{equation} 

It is very important that~\eqref{eq:8.4} holds for character sums over 
any interval of length $N$, not only for the initial segment $0<n \le N$. 
If $M=0$ the Riemann Hypothesis yields  
\begin{equation}\label{eq:8.7}
S_{\chi}(0,N) \ll N^{\frac 12}q^{\ve}\ ,
\end{equation} 
so~\eqref{eq:8.4} holds with $\delta = 1/2$, but this special case is 
not sufficient for our applications. We hope that Conjecture $C_n$ will 
be established in the not too distant future by advancing the tools 
of analytic number theory so our main theorem will become unconditional 
for fields of any degree $n \ge 3$. Note that by Burgess' 
result~\eqref{eq:8.2}, we narrowly missed Conjecture $C_4$. 

Recall, we actually needed these bounds for character sums over 
an arithmetic progression. Fortunately,
there is an immediate consequence of Conjecture $C_n$ for such sums:
\begin{equation}\label{eq:8.8}
S_{\chi}(M,N;k,\ell) = 
\sum_{\substack{M< n \le M+N\\ n\equiv \ell\, (\mod k)}} \chi(n)\ .
\end{equation} 
This can be written as
$$
\sum_{\frac{M-\ell}{k}< m \le \frac{M +N-\ell}{k}} \chi(km+\ell)\ .
$$ 
Now, suppose $(k,q)=1$. Take $\bar k$ with $\bar k k \equiv 1\, (\mod q)$. 
Then, we have
$$
\bar\chi (k)S_{\chi}(M,N;k,\ell) = 
\sum_{\frac{M-\ell}{k}< m \le \frac{M +N-\ell}{k}} \chi(m+\ell\bar k)
=S_{\chi}(M', N/k) 
$$ 
with $M' = (M-\ell)k^{-1} + \ell \bar k$. Since the bound~\eqref{eq:8.4} 
does not depend on $M$, it remains valid for $S_{\chi}(M,N;k,\ell)$ 
provided $\chi \, (\mod q)$ is not principal and $(k,q)=1$. 

\section{Relating Spins and Selmer Groups}\label{sec 11}

Suppose $E$ is an elliptic curve over $\qb$, and let $K = \qb(E[2])$, 
the number field 
generated by the points of order $2$ in $E(\bar\qb)$.  We assume that 
$K$ is a cyclic cubic extension 
of $\qb$ and that all totally positive units of $K$ are squares, 
so $K$ satisfies the 
hypotheses of the earlier sections of this paper.

Let 
$$
y^2 = f(x)
$$
be a Weierstrass model of $E$.  Then $f$ is irreducible over $\qb$, 
and we can identify $K$ with 
$\qb[T]/f(T)$.  For every place $v$ of $\qb$, we let 
$K_v = K \otimes \qb_v = \qb_v[T]/f(T)$.
If $d \in \qb^\times$, then the quadratic twist $E^{(d)}$ of $E$ 
by $d$ is the elliptic curve $dy^2 = f(x)$. 

See for example \cite[\S2]{MR} for the definition of the $2$-Selmer group 
$\mathrm{Sel}_2(E)$, or else take the description \eqref{selmerdef} 
below as the definition.  The $2$-Selmer group is an $\mathbb{F}_2$-vector 
subspace of $H^1(\qb,E[2])$, sitting in an exact sequence 
\newfont{\cyrr}{wncyr10}
\def\Sh{\mbox{\cyrr Sh}}
$$
0 \to E(\qb)/2E(\qb) \to \mathrm{Sel}_2(E) \to \Sh(E)[2]  \to 0
$$
where $\Sh(E)$ is the Shafarevich-Tate group of $E/\qb$.
Let 
$$
\Sigma = \{\text{primes $\ell$} : \text{$E$ has bad reduction at $\ell$, 
and $\ell$ is unramified in $K$}\} 
   \cup \{2\}.
$$
The main result of this section is the following.

\begin{theo}
\label{thm11.1}
Suppose $p$ is a prime that splits completely in $\qb(E[4])$, and let 
$\fp$ be a prime of $K$ above $p$.  
Suppose further that $\fp$ has a totally positive generator $\pi$ such 
that $\pi$ is a square  
in $K_\ell$ for every prime $\ell \in \Sigma$.
Then viewing $\mathrm{Sel}_2(E), \mathrm{Sel}_2(E^{(p)}) 
\subset H^1(\qb,E[2])$, we have 
$\mathrm{Sel}_2(E) \subset \mathrm{Sel}_2(E^{(p)})$ and 
$$
\dim_{\mathbb{F}_2}\mathrm{Sel}_2(E^{(p)}) = 
\begin{cases}
\dim_{\mathbb{F}_2}\mathrm{Sel}_2(E) + 2 & \text{if $spin(\fp) = +1$}, \\
\dim_{\mathbb{F}_2}\mathrm{Sel}_2(E) & \text{if $spin(\fp) = -1$}.
\end{cases}
$$
\end{theo}

To prove Theorem \ref{thm11.1}, we will  use the description of the 
$2$-Selmer group given in 
\cite{brumerkramer}.  Namely, let 
$$
(K^\times/(K^\times)^2)^{N=1} := \{\alpha \in K^\times/(K^\times)^2 : 
N_{K/\qb}(\alpha) \in (\qb^\times)^2\}
$$
and similarly for $(K_v^\times/(K_v^\times)^2)^{N=1}$, for every 
place $v$ of $K$.
For every $v$ there is a commutative diagram
$$
\xymatrix@C=40pt{
E(\qb)/2E(\qb) \ar@{^(->}^{\lambda_{E/\qb}}[r]\ar[d] 
& (K^\times/(K^\times)^2)^{N=1} \ar^{\iota_v}[d] \\
E(\qb_v)/2E(\qb_v) \ar@{^(->}^-{\lambda_{E/\qb_v}}[r] 
& (K_v^\times/(K_v^\times)^2)^{N=1}.
}
$$
The injection $\lambda_{E/\qb}$ is defined for $P \in E(\qb) - E(\qb)[2]$ 
by $\lambda_E(P) := x(P)-T$, 
where $x(P)$ denotes the $x$-coordinate of $P$ and we identify 
$K = \qb[T]/f(T)$, 
and $\lambda_{E/\qb_v}$ is defined similarly.
By \cite[\S2]{brumerkramer} there is a natural identification 
\begin{equation}
\label{selmerdef}
\mathrm{Sel}_2(E) = \{\alpha \in (K^\times/(K^\times)^2)^{N=1} : 
   \text{$\iota_v(\alpha) \in \mathrm{image}(\lambda_{E/\qb_v})$ 
for every $v$}\}.
\end{equation}

Fix a generator $\sigma$ of $\mathrm{Gal}(K/\qb)$.  If $\fp$ is a 
prime of degree one of $K$, 
lying above the rational prime $p$, we identify $K_p$ 
with $\qb_p\times \qb_p \times \qb_p$.
The map $\iota_p$ is then given by  
\begin{equation}
\label{iota}
\iota_p(\alpha) = (\alpha,\alpha^\sigma,\alpha^{\sigma^2}) 
   \in (K_\fp^\times/(K_\fp^\times)^2)^3 
   = (\qb_p^\times/(\qb_p^\times)^2)^3.
\end{equation}

\begin{lem}\label{lem 11.2}
Suppose  $p$ is an odd prime that splits completely in $\qb(E[4])$, 
and $E$ has good reduction at $p$.  
With the identification above, we have
\begin{align}
\label{l11.2.1}
&\text{$\mathrm{image}(\lambda_{E/\qb_p})$ is generated by $(u,u,1)$ 
and $(u,1,u)$}, \\
\intertext{where $u \in \zb_p^\times$ is a nonsquare, and}
\label{l11.2.2}
&\text{$\mathrm{image}(\lambda_{E^{(p)}/\qb_p})$ is generated 
by $(p,p,1)$ and $(p,1,p)$}.
\end{align}
\end{lem}

\begin{proof}
Assertion \eqref{l11.2.1} is \cite[Corollary 3.3]{brumerkramer}.

Since $p$ splits completely in $\qb(E[4])$, $\mathrm{Gal}(\bar{\qb}_p/\qb_p)$ 
acts trivially on 
$E[4]$, and hence $\mathrm{Gal}(\bar{\qb}_p/\qb_p)$ acts on $E^{(p)}[4]$ 
by the quadratic 
character of $\qb_p(\sqrt{p})/\qb_p$.  
It follows that the natural map 
$E^{(p)}(\qb_p)[2] \to E^{(p)}(\qb_p)/2E^{(p)}(\qb_p)$ is an isomorphism  
and the natural map
$$
E^{(p)}(\qb_p)/2E^{(p)}(\qb_p) \to E^{(p)}(\qb_p(\sqrt{p}))/2E^{(p)}
(\qb_p(\sqrt{p}))
$$
is identically zero.
Therefore, if $\lambda_{E^{(p)}/\qb_p}(P) = (t_1,t_2,t_3) 
\in (\qb_p^\times/(\qb_p^\times)^2)^3$, then 
$t_i$ is a square in $\qb_p(\sqrt{p})$ for every $i$.  
Since $|E^{(p)}(\qb_p)/2E^{(p)}(\qb_p)| = 4$, \eqref{l11.2.2} follows.
\end{proof}

\begin{proof}[Proof of Theorem \ref{thm11.1}]
Suppose $p$ satisfies the hypotheses of Theorem \ref{thm11.1}.
Then $p = N_{K/\qb}\pi$ is a square in $\qb_v$ if 
$v \in \Sigma \cup \{\infty\}$. 
Hence for such $v$, $E$ is isomorphic to $E^{(p)}$ over $\qb_v$, so
$\mathrm{image}(\lambda_{E/\qb_v}) = \mathrm{image}(\lambda_{E^{(p)}/\qb_v})$.
If $\ell \ne 2$ ramifies in $K$, then $E(\qb_\ell)[2] 
= E^{(p)}(\qb_\ell)[2] = 0$, 
so $E(\qb_\ell)/2E(\qb_\ell) = E^{(p)}(\qb_\ell)/2E^{(p)}(\qb_\ell) = 0$ and 
$\mathrm{image}(\lambda_{E/\qb_\ell}) 
= \mathrm{image}(\lambda_{E^{(p)}/\qb_\ell}) = 0$.  
If $\ell \nmid 2p$ is a prime where $E$ has good reduction, 
unramified in $K$, then 
\cite[Lemma 4.1]{cassels8} shows that 
$\mathrm{image}(\lambda_{E/\qb_\ell}) = \mathrm{image}
(\lambda_{E^{(p)}/\qb_\ell})$.  Therefore 
$\mathrm{image}(\lambda_{E/\qb_v}) = \mathrm{image}
(\lambda_{E^{(p)}/\qb_v})$ for every $v \ne p$.

Let 
\begin{align*}
S^{(p)} &= \{\alpha \in (K^\times/(K^\times)^2)^{N=1} : 
   \text{$\iota_v(\alpha) \in \mathrm{image}(\lambda_{E/\qb_v})$ 
for every $v \ne p$}\}, \\
S_{(p)} &= \{\alpha \in S^{(p)} : \iota_p(\alpha) = 1\}.
\end{align*}
Then 
\begin{align*}
\mathrm{Sel}_2(E) 
   &= \{\alpha \in S^{(p)} : \iota_p(\alpha) \in \mathrm{image}
(\lambda_{E/\qb_p})\},\\
\mathrm{Sel}_2(E^{(p)}) 
   &= \{\alpha \in S^{(p)} : \iota_p(\alpha) \in \mathrm{image}
(\lambda_{E^{(p)}/\qb_p})\}.
\end{align*}
By Lemma~\ref{lem 11.2}, $\mathrm{image}(\lambda_{E/\qb_p}) \cap 
\mathrm{image}(\lambda_{E^{(p)}/\qb_p}) = \{1\}$, 
so 
$$
\mathrm{Sel}_2(E) \cap \mathrm{Sel}_2(E^{(p)}) = S_{(p)}.
$$  
Global duality (see for example \cite[Lemma 3.2]{MR}) shows that
\begin{equation}
\label{duality}
\dim_{\mathbb{F}_2} \iota_p(S^{(p)}) = \dim_{\mathbb{F}_2} 
(S^{(p)}/S_{(p)}) = 2.
\end{equation}

Let $\pi$ be the totally positive generator of a prime above $p$, 
as in the statement of Theorem 
\ref{thm11.1}, and let $\alpha = \pi\pi^\sigma$.  
Then $N_{K/\qb}\alpha = p^2$.  
Since $\pi$ is totally positive, $\alpha \in (K_\infty^\times)^2$.  
By assumption, 
$\alpha \in (K_\ell^\times)^2$ if $\ell = 2$ or if $\ell$ is a prime 
of bad reduction, unramified 
in $K$.  If $\ell$ is odd and $\ell$ ramifies in $K$, then 
$\alpha \equiv \pi^2$ modulo the prime of 
$K$ above $\ell$, so $\alpha \in (K_\ell^\times)^2$ in this case 
as well.  Finally, if $\ell \ne p$ 
then $\alpha$ is a unit at $\ell$, so if $\ell$ is a prime of good 
reduction then (again using 
\cite[Lemma 4.1]{cassels8}) we have $\iota_\ell(\alpha) 
\in \mathrm{image}(\lambda_{E/\qb_\ell})$.  
Therefore $\alpha \in S^{(p)}$.

Since $\iota_p(\alpha)$ and $\iota_p(\alpha^\sigma)$ are distinct 
and nonzero in $K_p^\times/(K_p^\times)^2$, 
we see by \eqref{duality} that $\alpha, \alpha^\sigma$ generate 
$S^{(p)}/S_{(p)}$.
By \eqref{l11.2.1}, $\iota_p(\alpha)$ and $\iota_p(\alpha^\sigma)$ generate 
a subgroup visibly disjoint from $\mathrm{image}(\lambda_{E/\qb_p})$.  Thus
$$
\iota_p(S^{(p)}) \cap \mathrm{image}(\lambda_{E/\qb_p}) = \{1\}
$$
so $\mathrm{Sel}_2(E) = S_{(p)}$.

By \eqref{iota} we have
$$
\iota_p(\alpha) = (\pi\pi^\sigma,\pi^\sigma\pi^{\sigma^2},\pi\pi^{\sigma^2}) 
   = (p,1,p)\cdot(1/\pi^{\sigma^2},\pi^\sigma\pi^{\sigma^2},1/\pi^\sigma) 
\in K_p^\times/(K_p^\times)^2\ .
$$
Since $\pi$ is a square in $K_2^\times$, Lemma~\ref{lem 9.1} in 
the next section shows that 
$\pi^\sigma, \pi^{\sigma^2} \in (K_\fp^\times)^2$ if and only if 
$spin(\fp) = +1$. 
Thus, by \eqref{l11.2.2}  
$$
\iota_p(\alpha) \in \mathrm{image}(\lambda_{E^{(p)}/\qb_p}) 
\iff spin(\fp) = +1\ .
$$
The same holds for $\alpha^\sigma$ and $\alpha\alpha^\sigma$, so we have 
$$
\mathrm{Sel}_2(E^{(p)}) = 
\begin{cases}
S^{(p)} & \text{if $spin(\fp) = +1$}, \\
S_{(p)} & \text{if $spin(\fp) = -1$}.
\end{cases}
$$
Now the theorem follows from \eqref{duality}.
\end{proof}

\noindent{\bf Example.}
Let $E$ be the elliptic curve $y^2 = x^3+x^2-16x-29$ of conductor 
$784$, as in the introduction.
Then $\dim_{\mathbb{F}_2}\mathrm{Sel}_2(E) = 1$,  
$K = \qb(E[2])$ is the real subfield of the field of $7$-th roots of 
unity, and every 
totally positive unit of $K$ is a square.

The only prime of bad reduction for $E$ that is unramified in $K$ is $2$, 
and $\qb(E[4])$ is contained 
in the ray class field of $K$ modulo $8\infty_1\infty_2\infty_3$.  Hence  
if $p$ is a rational prime that splits completely in $K$, and a prime 
$\fp$ above $p$ 
has a totally positive generator 
congruent to $1$ modulo $8$, then $p$ splits completely in $\qb(E[4])$ 
and 
Theorem \ref{thm11.1} applies to show that
$$
\dim_{\mathbb{F}_2}\mathrm{Sel}_2(E^{(p)}) = 
\begin{cases}
3 & \text{if $spin(\fp) = +1$}, \\
1 & \text{if $spin(\fp) = -1$}.
\end{cases}
$$

\section{More than One Spin}\label{sec 9} 

There is plenty of scope for further work on the subject of this paper. 
We can for example ask what happens if one or more of the 
various restrictions are dropped from those we placed on the number fields 
being considered: 
``totally real, Galois, cyclic, units achieving all signs''. 

Even keeping these assumptions in place, there remain interesting problems.
A particularly natural one is the question of the joint distribution 
of the $n$ different spins attached to the primes by the various 
automorphisms of $G= Gal (K/\qb)$. An instantaneous guess might be that 
the distributions 
are independent of each other but, even in the simplest case of cyclic $G$,  
there are relations amongst the spins which make the situation a little 
more complicated, and perhaps more interesting. 

Let $\sigma$ and $\tau$ be two elements of $G$ which, given our 
assumptions, commute. Let $\fa =(\alpha)$ where $\alpha$ is odd and totally 
positive. Recovering our notation from~\eqref{eq:1.1}, we have 
$$
spin(\sigma,\fa^{\tau}) 
= \Bigl(\frac{\alpha^{\tau}}{(\alpha)^{\tau\sigma}}\Bigr)
= \Bigl(\frac{\alpha^{\tau}}{(\alpha)^{\sigma\tau}}\Bigr)
= \Bigl(\frac{\alpha}{(\alpha)^{\sigma}}\Bigr) = spin(\sigma,\fa)\ , 
$$
a fact we have already been using.
Specializing to $\tau=\sigma^{-1}$ we find that
$$
spin(\sigma^{-1},\fa) = \Bigl(\frac{\alpha}{(\alpha)^{\sigma^{-1}}}\Bigr)
= \Bigl(\frac{\alpha^{\sigma}}{(\alpha)^{\sigma^{-1}\sigma}}\Bigr) 
= \Bigl(\frac{\alpha^{\sigma}}{\alpha}\Bigr)\ , 
$$
so that, by quadratic reciprocity: 
\begin{lem}\label{lem 9.1} 
We have 
$$
spin(\sigma, \frak a) 
= spin (\sigma^{-1}, \frak a) \mu_2(\alpha, \alpha^{\sigma})\ , 
$$ 
where 
$\mu_2(\alpha, \alpha^{\sigma})$   is a product of the local Hilbert 
symbols at primes $\frak p \mid 2$.  In particular, 
$\mu_2(\alpha, \alpha^{\sigma}) = 1$  if   $\alpha \equiv  1\,( \mod 4)$.
\end{lem}

If we now let $\sigma$ be a generator of $G$ then, for each 
$k,\, 1 \le k \le n-1$, 
there is an evident dependence between $spin(\sigma^k,\fa)$ 
and $spin(\sigma^{n-k},\fa)$. For fields of odd degree this shows that 
there is no non-identity  
spin whose distribution is independent of all of the others 
and the question arises as to whether these are the only dependencies 
among the spins. For the simplest possible example we can ask:

\noindent{\bf Problem}. 
For $K$ totally real, of degree $n\ge 5$ (odd or even), with all totally 
positive units being squares, and $G= Gal (K/\qb)$ being cyclic and 
generated by $\sigma$, are the distributions of $spin(\sigma, \fa)$
and $spin(\sigma^2, \fa)$ independent?

If we consider fields $K$ of even degree the (non-identity) 
spins also pair off and the same questions 
arise, apart from the middle spin, $k=n/2$. Here, we have an involution, 
say $\sigma =\sigma^{-1}$, and a somewhat different picture emerges. 
Recall that, in addition to our assumptions about the field, we have 
also been restricting our attention to generators of the Galois group 
so that, even though we are now talking about a single spin, 
we have not dealt with this situation. In Section~\ref{sec 10}  
using very different arguments, we give results  
which show that equidistribution holds for the sum of 
prime spins attached to such an involution.

\section{Prime Spins for an Involution}\label{sec 10}  

In this section we are going to require our Galois automorphism $\sigma$  
to be an involution (not the identity), rather than a generator of the 
group so, in addition 
to our previous restrictions, we ask that $K$ have even degree $n$ and 
we denote by $L$ the fixed field of $\sigma$ so that $K/L$ is a quadratic 
field extension. Thus, in particular, $K$ could be a real quadratic field 
so long as its fundamental unit has norm $-1$. 
We shall also make some further simplifications as 
we proceed. 

Our first main result gives a natural arithmetic characterization 
of the spin attached to such an involution. 
 
\begin{prop}\label{prop 10.1} 
Let $K/\qb$ be a totally real Galois extension of even degree $n\ge 2$ with 
${\cal U}^+  = {\cal U}^2$. 
Let $L$ be the subfield of $K$ fixed by the involution $\sigma$ and suppose 
that the discriminant $\fD$ of the quadratic extension $K/L$ is odd. 
Let $\alpha \in {\cal O}_K,$ with $(\alpha, \, \alpha^{\sigma})=1, \, 
\alpha \equiv 1\, (\mod 8), \, \alpha \succ 0$. Then, for $\fa = (\alpha)$, 
we have
\begin{equation}\label{eq:10.0}
spin (\sigma, \fa) = \Bigl( \frac{\beta}{\fD}\Bigr)_L 
\end{equation}
where
$$ 
\beta= {\tfrac 12} (\alpha + \alpha^{\sigma}) 
={\tfrac 12} T_{K/L}(\alpha)\ .
$$ 
\end{prop}
 Note that the notation in the quadratic symbol 
now reflects the field. 
This reduction of the spin to something so close to the field character 
imparts multiplicativity properties which render inoperable the method 
we have employed up to now
but open the possibility of using the theory of $L$-functions. As a result,  
we shall obtain the following theorem. 
 
\begin{theo}\label{theo 10.2}  
Let $K/\qb$ be a totally real Galois extension of even degree $n\ge 2$ with 
${\cal U}^+  = {\cal U}^2$. Let $\sigma \in Gal(K/\qb)$ be an involution 
and let the discriminant of the relative quadratic extension be odd.
Let $\fP$ run over the principal prime ideals of $K$ with
$$
\fP=(\alpha), \,\, \alpha \succ 0,\,\, \alpha \equiv 1\, (\mod 8)\ .
$$
Then, we have 
$$
\sum_{N\fP \le x} spin (\sigma, \fP) \ll x\exp (-C\sqrt{\log x})\ ,
$$
where the positive constant $C$ and the implied constant depend 
on the field $K$.
On assumption of the Riemann Hypothesis for Hecke $L$-functions, 
the above bound can be sharpened to
$$
\sum_{N\fP \le x} spin (\sigma, \fP) \ll x^{\frac 12} (\log x)^A\ , 
$$
where now the positive constant $A$ depends on the degree of the field 
$K$. 
\end{theo}
\noindent
Note that we now denote the prime ideals of $K$ by $\fP$ rather than by 
$\fp$ as before, because we reserve $\fp$ for the prime ideals of the subfield
$L$. 

The path to Proposition~\ref{prop 10.1} passes through a number of lemmas 
which lead to progressively simpler expressions for the spin. 
We begin with 
\begin{lem}\label{lem 10.3}
If $\fP \subset {\cal O}_K$ is an odd prime with $(\fP, \fP^{\sigma}) = 1$ 
then, for $x\in {\cal O}_L$ we have
\begin{equation}\label{eq:10.1}
\Bigl(\frac{x}{\fP}\Bigr)_K =\Bigl(\frac{x}{\fp}\Bigr)_L\, 
\quad {\rm with}\,\, \fp =\fP \fP^{\sigma}\ . 
\end{equation} 
\end{lem} 
\begin{proof}
This follows quickly on combining the Euler criteria for the two fields:
\begin{equation*}
\Bigl(\frac{x}{\fP}\Bigr)_K  \equiv x^{\frac 12 
(N_{K/\qb}(\fP) -1)}\, (\, \mod \fP)\ , 
\end{equation*}
\begin{equation*}
\Bigl(\frac{x}{\fp}\Bigr)_L  \equiv x^{\frac 12 
(N_{L/\qb}(\fp) -1)}\, (\, \mod \fp)\ ,
\end{equation*}
valid for $(x, \fp)=1$. Here, $N_{K/\qb}(\fP) = N_{L/\qb}(\fp)$. 
\end{proof}
\begin{cor}\label{cor 10.4} If $\alpha \in {\cal O}_K$ is odd and 
$(\alpha, \alpha^{\sigma})= 1$ then, for $x\in {\cal O}_L$, we have 
\begin{equation}\label{eq:10.2}
\Bigl(\frac{x}{\alpha}\Bigr)_K 
=\Bigl(\frac{x}{\alpha\alpha^{\sigma}}\Bigr)_L\ .
\end{equation} 
\end{cor}

Now, let $\alpha \in {\cal O}_K$ be odd, $(\alpha, \alpha^{\sigma})= 1$ 
and $\alpha\succ 0$. Then, the Corollary gives
$$
spin (\sigma, \alpha) = \Bigl(\frac{\alpha}{\alpha^{\sigma}}\Bigr)_K 
= \Bigl(\frac{\alpha + \alpha^{\sigma}}{\alpha^{\sigma}}\Bigr)_K 
= \Bigl(\frac{\alpha + \alpha^{\sigma}}{\alpha \alpha^{\sigma}}\Bigr)_L\ .
$$ 
Suppose that $\alpha \equiv 1\, (\mod 8)$. Then,
$$
\beta = \tfrac12(\alpha + \alpha^{\sigma}) \equiv 1\, (\mod 4)\ , \quad
\gamma  = \tfrac12(\alpha - \alpha^{\sigma}) \equiv 0\, (\mod 4)\ , 
$$
satisfy $(\beta,\gamma)=1, \, \beta \succ 0$ and 
$\alpha \alpha^{\sigma} =\beta^2 - \gamma^2 \equiv 1 \, (\mod 8)$. 
We have $\beta \in {\cal O}_L$ and $\gamma^2 \in {\cal O}_L$. Hence, the spin 
simplifies to 
\begin{equation}\label{eq:10.4}
spin (\sigma, \alpha) = 
\Bigl(\frac{2}{\alpha \alpha^{\sigma}}\Bigr)_L
\Bigl(\frac{\beta}{\alpha \alpha^{\sigma}}\Bigr)_L
 = \Bigl(\frac{2}{\alpha \alpha^{\sigma}}\Bigr)_L
\Bigl(\frac{\alpha \alpha^{\sigma}}{\beta}\Bigr)_L
=\Bigl(\frac{-\gamma^2}{\beta}\Bigr)_L\ ,
\end{equation}
since $\alpha \equiv 1\, (\mod 8)$ implies 
$\alpha^{\sigma} \equiv 1\, (\mod 8)$, which implies
$ (\frac{2}{\alpha \alpha^{\sigma}})_L= 1$. 

The final lemma in this chain is:
\begin{lem}\label{lem 10.5} If $\beta \equiv 1\, (\mod 4)$ 
and $\delta$ odd, are in ${\cal O}_L$, then 
$$
\Bigl(\frac{\beta, \, \delta}{\fp}\Bigr) 
= 1\quad {\rm for\,\, every} \quad \fp \mid 2\ . 
$$
\end{lem}
\begin{proof} First we consider the special case 
$\delta \equiv 1\, (\mod 2)$. Here, we use 
the identity $1 + 4x =(1+2y)^2 -4\delta y^2$ with
$x = y + (1-\delta )y^2$.  
Given $x \in {\cal O}_L$ the latter equation is solvable for 
$y \in L_{\fp}$ by Hensel's Lemma. This yields the result. 

In general, if $\delta$ is odd, there exists an odd positive integer $r$ 
such that $\delta^r \equiv 1\, (\mod 2)$. This follows from Fermat's Little 
Theorem with $r=N\fp -1$ for any $\fp \mid 2$. Hence, the result follows from 
the special case by the multiplicativity of the Hilbert symbol.  
\end{proof} 

Now, consider the submodule ${\cal M} = {\cal O}_L + \frac{(1+\alpha)}{2}
{\cal O}_L $ 
of ${\cal O}_K$. Since 
$$
{\rm det}
\begin{pmatrix}
1 & (1+\alpha)/2\\ 1 & (1+\alpha^{\sigma})/2\end{pmatrix} 
= - \, \gamma\ , 
$$
the discriminant of ${\cal M}$ is the principal ideal $(\gamma^2)$ of 
${\cal O}_L$. Hence, $(\gamma^2)=\fa^2\fD$ where $\fa$ is an ideal of 
${\cal O}_L$ such that ${\cal O}_K/{\cal M} \cong {\cal O}_L/\fa$ and 
$\fD$ is the discriminant of the relative quadratic extension $K/L$. 
We have $(\beta,\gamma^2)=1$ so $(\beta, \fa^2 \fD)=1$. Choose an ideal 
$\fb$ of ${\cal O}_L$ with $\fb \sim \fa$ and $(\fb, 2\beta\fa)=1$ 
and choose $\mu\in L$ 
such that $(\mu) = \fa^{-1}\fb$. Put $\gamma_1 = \gamma\mu \in K$. Then, 
$\gamma_1^2 = \gamma^2\mu^2 \in L$ and, in fact,
\begin{equation}\label{eq:10.5}
(\gamma_1^2) = (\gamma^2)\,(\mu)^2 =\fa^2\fD\,(\fa^{-1}\fb)^2 =\fb^2\fD\ ,
\end{equation} 
so $\gamma_1^2 \in {\cal O}_L$. Since $(\beta, \fa\fb )=1$ we have 
$$
\Bigl(\frac{-\gamma_1^2}{\beta}\Bigr)_L = 
\Bigl(\frac{-\gamma^2\mu^2}{\beta}\Bigr)_L =
\Bigl(\frac{-\gamma^2}{\beta}\Bigr)_L\ . 
$$
By~\eqref{eq:10.5}, $(\gamma_1^2,2\beta)=1$ and 
$$
\Bigl(\frac{-\gamma_1^2}{\beta}\Bigr)_L = 
\Bigl(\frac{\beta}{\gamma_1^2}\Bigr)_L \prod_{\fp\mid 2\infty}
\Bigl(\frac{\beta, -\gamma_1^2}{\fp}\Bigr)
= \Bigl(\frac{\beta}{\gamma_1^2}\Bigr)_L 
=\Bigl(\frac{\beta}{\fD}\Bigr)_L\ ,
$$ 
by the reciprocity law and Lemma~\ref{lem 10.5}. 
Now, by~\eqref{eq:10.4}, this completes the proof of 
Proposition~\ref{prop 10.1}. 

We now turn to the proof of Theorem~\ref{theo 10.2}.
\begin{proof}
Let $\fP$ run over the principal prime ideals of $K$ with
\begin{equation}\label{eq:10.8}
\fP = (\alpha), \,\,\, \alpha\succ 0, \,\,\, \alpha \equiv 1\, (\mod 8)\ .
\end{equation}
Denote by $S(x)$ the number of such prime ideals with $N\fP \le x$ 
and by $S(x; \fD, \delta)$ the number of these in the residue 
class $\alpha \equiv \delta \, (\mod \fD)$. 
By the Prime Ideal Theorem we have (recall that $\fD$ is odd) 
$$
S(x; \fD, \delta) \sim \frac{S(x)}{\varphi_K(\fD)}, 
\quad {\rm if} \,\,\, (\delta, \fD) = 1\ , 
$$ 
where $\varphi_K(\fD)$ is the number of classes 
$\delta \, (\mod \fD)$ in ${\cal O}_K$ with $(\delta, \fD) = 1$.
More precisely, the error term satisfies 
\begin{equation}\label{eq:10.9}
E(x; \fD, \delta) = S(x; \fD, \delta)\, 
-\, \frac{S(x)}{\varphi_K(\fD)} \ll x\exp\bigl(- C \sqrt {\log x} \bigl)
\end{equation}
unconditionally, and 
\begin{equation}\label{eq:10.10}
E(x; \fD, \delta) \ll x^{1/2}(\log x)^A
\end{equation} 
subject to the Riemann Hypothesis for the relevant Hecke $L$-functions. 
Here $C$ and $A$ are positive constants depending on the field $K$, 
as do the implied constants.  

Note that
\begin{equation}\label{eq:10.11}
\sum_{\substack{\delta\, (\mod\fD)\\ (\delta, \fD) = 1}}
\Bigl(\frac{\delta + \delta^{\sigma}}{\fD}\Bigr)_L = 0\ .
\end{equation} 
To see this we change the variable $\delta$ to $\delta\eta$ with 
$\eta\in {\cal O}_L, \, (\eta, \fD)=1$. We find the sum is equal to 
$(\eta/\fD)_L$ times itself. Choosing $\eta$ such that 
$\chi_{\fD}(\eta) = (\eta/\fD)_L = -1$, we obtain~\eqref{eq:10.11}. 

Now, put 
\begin{equation}\label{eq:10.13}
S^{\sigma} (x) =\sum_{N\fP \le x} spin (\sigma, \fP)\ ,
\end{equation} 
where $\fP$ runs through the prime ideals~\eqref{eq:10.8}. 
By the formula~\eqref{eq:10.0} we get
$$
S^{\sigma} (x) = \Bigl(\frac{2}{\fD}\Bigr)_L
\sum_{\substack{\delta\, (\mod\fD)\\ (\delta, \fD) = 1}}
\Bigl(\frac{\delta + \delta^{\sigma}}{\fD}\Bigr)_L\, S(x; \fD, \delta) 
+ O(1)\ ,
$$
where the error term $O(1)$ takes into account the contribution of 
the prime ideals dividing $\fD$. Finally, applying the Prime Ideal 
Theorem in the form~\eqref{eq:10.9} or the Riemann Hypothesis in the 
form~\eqref{eq:10.10}, we find by~\eqref{eq:10.11} that the main terms 
cancel out so we are left with the result claimed in the theorem.
\end{proof}

%
%
%


\begin{thebibliography}{xxxx}

\bibitem[BK]{brumerkramer} 
A.\ Brumer, K.\ Kramer, The rank of elliptic curves,  
Duke Math. J. 44 (1977) 715-Ð743. 

\bibitem[Bu]{Bu} D. A. Burgess, On character sums and {$L$}-series {II}, 
Proc. London Math. Soc., 13, (1963), 524--536. 

\bibitem[Ca1]{cassels8}
J.W.S.\ Cassels, Arithmetic of curves of genus 1, VIII,  
J. Reine Angew.\ Math. 217 (1965) 180-Ð199.


\bibitem[Ca2]{Ca}  J. W. S. Cassels, 
Global fields, Algebraic Number Theory (Proc. Instructional Conf., 
Brighton, 1965) pp. 42--84 Thompson, (Washington) 1967. 

\bibitem[DFI]{DFI} W. Duke, J. B. Friedlander and H. Iwaniec, 
Bilinear forms with Kloosterman fractions, 
Invent. Math. 128 (1997), 23--43. 

\bibitem[FI]{FI}  J. B. Friedlander and H. Iwaniec, The polynomial 
$X^2+Y^4$ captures its primes, Ann. Math. 148 (1998), 945--1040. 

\bibitem[He]{He} E. Hecke, Lectures on the Theory of Algebraic Numbers, 
Grad. Texts in Math. 77, Springer (New York) 1981.

\bibitem[IK]{IK} H. Iwaniec and E. Kowalski, Analytic Number Theory, 
Colloq. Pub. 53, Amer. Math. Soc. (Providence) 2004.

\bibitem[La]{La} S. Lang, Algebraic Number Theory, 2nd ed., 
Grad. Texts in Math. 110, Springer (New York) 1994. 


\bibitem[Le]{Le} E. Lehmer, Connection between Gaussian periods and cyclic
units, Math. Comp. 50 (1988) 535--541.

\bibitem[MR]{MR} B. Mazur, B. and K. Rubin, Ranks of twists of 
elliptic curves and Hilbert's tenth problem, Invent. Math. 181 (2010) 
541--575. 


\bibitem[Ne]{Ne} J. Neukirch, Algebraic Number Theory.
 Grund. Math. Wiss. 322. Springer (Berlin) 1999. 

\bibitem[Sha]{Sha} D. Shanks, The simplest cubic fields, Math. Comp. 28 
(1974) 1137-1152.

\bibitem[Shi]{Shi} T. Shintani, On evaluation of zeta functions of totally 
real algebraic number fields at non-positive integers, 
J. Fac. Sci. Univ. Tokyo Sect. IA Math. 23 (1976) 393--417. 

\end{thebibliography}
\end{document}